\documentclass[preprint,12pt]{elsarticle}
\usepackage{amsfonts}
\usepackage{amsmath}
\usepackage{amssymb}
\usepackage{amsthm}
\usepackage{setspace}
\usepackage[margin=1.0in]{geometry}
\usepackage{tikz}
\usepgflibrary{patterns}

\usepackage[draft]{todonotes}   

\DeclareMathOperator{\sign}{sign}

\usepackage{footnote}
\makesavenoteenv{tabular}
\makesavenoteenv{table}

\usepackage{hyperref}

\newtheorem{theorem}{Theorem}
\newtheorem{definition}{Definition}

\newtheorem{lemma}{Lemma}
\newtheorem{remark}{Remark}
\newtheorem{proposition}{Proposition}
\newtheorem*{proof*}{Proof}

\protected\edef\mathbb{%
  \unexpanded\expandafter\expandafter\expandafter{%
    \csname mathbb \endcsname
  }%
}

\usepackage{amsthm}

\date{\today}

\journal{Journal of Computational and Applied Mathematics}

\begin{document}

\begin{frontmatter}
\title{An efficient approach for solving stiff nonlinear boundary value problems}
\author{Volodymyr L. Makarov}
\ead{makarov@imath.kiev.ua}
\author{Denys V. Dragunov}
\ead{dragunovdenis@gmail.com, dragunovdenis@imath.kiev.ua}
\address{Department of Numerical Mathematics \\ Institute of Mathematics \\ National Academy of Sciences of Ukraine\\ 01004 Ukraine, Kiev-4, 3, Tereschenkivska st.}

\begin{abstract}
A new method for solving stiff two-point boundary value problems is described and compared to other known approaches using the Troesch's problem as a test example. The method is based on the general idea of alternate approximation of either the unknown function or its inverse and has a genuine "immunity" towards numerical difficulties invoked by the rapid variation (stiffness) of the unknown solution. A c++ implementation of the proposed method is available at \url{https://github.com/imathsoft/MathSoftDevelopment}.
\end{abstract}


\end{frontmatter}

\section{Introduction}
In the present paper we consider a nonlinear boundary value problem (BVP)
\begin{eqnarray}
\frac{d^{2} u(x)}{d x^{2}} = N\left(u(x), x\right)u(x), \quad x\in\left[a, b\right], \; N(u,x)\in C^{2}(\mathbb{R}\times [a, b]), \label{first_eq} \\
u(a) = u_{a} \in \mathbb{R}, \quad u(b) = u_{b} \in \mathbb{R}, \label{boundary_cond}
\end{eqnarray}
which arises in many areas of physics and mathematics. Although, there is a huge variety of known methods for solving problems of type \eqref{first_eq}, \eqref{boundary_cond} (see, for example \cite{Allgower_mesh_refinements}, \cite{Chang20103303}, \cite{Ha20011411}, \cite{HERMANN1993103} and the references therein),
almost none of them fill comfortable when the problem turns out to be stiff.

As it was pointed out in \cite{Cash_Review_Paper_Num_meth_for_Stiff_IVP},
a good mathematical definition of the concept of stiffness does not exist. The famous definition given in \cite{Hairer_Wanner_Stiff} says that "stiff equations
are problems for which explicit methods don't work", which, unfortunately, is not very constructive. According to \cite{Fifty_years_of_stiffness}, there is at least 6 different definitions of
stiff problems which possess different levels of formality and are accepted by different schools of mathematics. The authors of \cite{Fifty_years_of_stiffness} came up with
their own definition of "stiffness", based on the concept of \textit{stiffness ratio}, which encompasses all the known definitions.

In the present paper we confine ourselves to consider only a subclass of
stiff boundary value problems \eqref{first_eq}, \eqref{boundary_cond} whose stiffness is originated from the fact that the exact solution $u(x)$
possesses narrow intervals of rapid variation, known as the  \textit{boundary layers}. Such a behavior it typical for \textit{singularly perturbed problems},
which are an important subclass of stiff problems (see,  \cite{Attili_Basem_2005}, \cite{Attili_Basem_2011}, \cite{Flaherty_O'Malley_stiff_differential_equations},  \cite{Flaherty_OMalley_1984}, \cite{Gaiduk_Slushaenko_large_L}, \cite{Kreiss_Nichols_Brown},  \cite{Fifty_years_of_stiffness}). The rapid variation
is equivalent to having $|u^{\prime}(x)| \gg 1$ on some subset of $[a, b].$ And it is the need to approximate the solution on this subset that makes the
problem numerically difficult and unstable, i.e. stiff. Now to approximate the solution on the subset of $[a, b]$ where $|u^{\prime}(x)|$ is comparatively small
is much easier from the numerical point of view. To be more specific, let us consider a set $\chi_{u} \in [a, b]$ defined in the following way:
\begin{equation}\label{Label_set_chi}
  \chi_{u} = \{x\in [a, b] : | u^{\prime}(x)| \geq 1\}.
\end{equation}
It is easy to see that, defined in such a way, set $\chi_{u}$ consists of a finite or infinite number of distinctive closed intervals $\bar{\iota}_{i}$. Some
of the  intervals
$\bar{\iota}_{i}$ might be those of rapid variation for the solution $u(x).$ At the same time, by the definition of $\chi_{u}$ \eqref{Label_set_chi}, solution $u(x)$
is strictly monotonic on each interval $\bar{\iota}_{i},$ which means that we can consider the inverse function $x_{\bar{\iota}_{i}}(\cdot) = u^{-1}(\cdot)$ defined on the closed interval
$u(\bar{\iota}_{i}) \in u([a, b]).$ There are two remarkable things about the function $x_{\bar{\iota}_{i}}(\cdot):$
\begin{enumerate}
  \item $|x^{\prime}_{\bar{\iota}_{i}}(u)| \leq 1, \forall u \in u(\bar{\iota}_{i}),$ which means that the initial BVP stated in terms of "inverse solution" $x^{\prime}_{\bar{\iota}_{i}}(u)$ is not stiff on $u(\bar{\iota}_{i});$
  \item having function $x_{\bar{\iota}_{i}}(u)$ approximated on a discrete set of points from $u(\bar{\iota}_{i})$ we automatically get function $u(x)$ approximated on some
  discrete set of points from $\bar{\iota}_{i}.$
\end{enumerate}
The two observations give us the key insight on how to deal with the subclass of stiff problems defined above. It is the \textit{divide and conquer} principle:
on the subintervals where solution $u(x)$ is well behaved (showing rather moderate variation) we solve the given problem \eqref{first_eq}, \eqref{boundary_cond},
whereas on the subintervals $\bar{\iota}_{i},$ where $u(x)$ varies rapidly (and the initial problem is stiff), we solve the corresponding problem for
the inverse solution $x_{\bar{\iota}_{i}}(u).$  Of course, this becomes feasible from the practical point of view only if there is a finite number of
subintervals $\bar{\iota}_{i},$ which becomes our assumption from now on.

Speaking about the known methods for solving BVPs, it is impossible not to mention the \textit{simple shooting method} (SSM) and the
\textit{multiple shooting method} (MSM) \cite[Section 7.3]{stoer2002introduction} which are two the most simple and reliable techniques to deal with
boundary value problems of type \eqref{first_eq}, \eqref{boundary_cond}. By calling them \textit{techniques} and not just \textit{methods} we would like
to emphasize that the basic idea behind them is very broad and can be used in many different modifications, which, in turn, might be called \textit{methods}.
Since definitions of both SSM and MSM essentially relay on using methods for solving \textit{initial value problems} (IVP), one of the ways to come up with a new
modification consists in using a different IVP solver. Below we adapt (modify)
the SSM and MSM by using a specific approach for numerical solution of
IVP's which is based on the idea of alternate approximation of either straight $u(x)$
or inverse $x(u)$ solutions of equation \eqref{first_eq} and has a genuine "immunity" towards numerical difficulties
invoked by the rapid variation (stiffness) of the solution in question.

The main focus of the paper is not only to present a general idea about how to treat some subclass of stiff boundary value problems in an efficient way, but also
to describe and examine a possible particular implementation of the idea, hereinafter referred to as {\it Straight-Inverse method} (or, simply, SI-method).
With this in mind, we actively exploit one of the most famous examples of stiff BVPs, known as the Troesch's problem:
\begin{eqnarray}
\frac{d^{2} u(x)}{d x^{2}} = \lambda \sinh\left(\lambda u(x)\right), \quad x\in\left[0, 1\right] \label{troesch_eq} \\
u(0) = 0, \quad u(1) = 1, \label{troesch_bound_cond}
\end{eqnarray}
which is a partial case of problem \eqref{first_eq}, \eqref{boundary_cond} with $N(u(x), x) \equiv \lambda\sinh\left(\lambda u(x)\right)/u(x),$ $a = u_{a}=0,$ $b = u_{b} = 1.$ In addition to its application in physics of plasma, the Troesch's problem, has drown a lot of interest to itself as a test case for methods of solving unstable
two-point boundary value problems because of its difficulties \cite{Snyman_Troesch_problem_solution}. A vast amount of numerical data available for the problem (see \cite{Snyman_Troesch_problem_solution}, \cite{Fifty_years_of_stiffness}, \cite{Chang20103043},
 \cite{JONES1973429}, \cite{Troesch1976279}, \cite{Gen_sol_of_TP} and the references therein) allowed us to perform broad analysis of the SI-method and compare it to many other methods for solving two-point BVPs. The comparison confirms excellent characteristics of the method in terms of both accuracy (numerical stability) and performance. Results of multiple numerical tests with problems other than \eqref{troesch_eq}, \eqref{troesch_bound_cond} (among them those with $N^{\prime}_{x}(u,x) \not\equiv 0$ and with the solution $u(x)$ oscillating on $[a, b]$), which are not included in the present paper, show remarkable adaptivity potential of the SI-method, and do support the conclusions obtained on the Troesch's test problem.

At this point, we would like to notice that, in general case, there is no guarantee that the BVP \eqref{first_eq}, \eqref{boundary_cond} is solvable, i.e. has a solution. From \cite[Theorem 7.25]{Kelley_Peterson_2010} it follows, however, that, under the conditions imposed on the nonlinearity of equation \eqref{first_eq}, the problem can have at most one solution, that is, the uniqueness is granted. The question of existence is kept out of the scope of the current paper, as well as the error analysis for the SI-method applied to BVP \eqref{first_eq}, \eqref{boundary_cond}. We leave both issues for the future publications. The main theoretical result of the paper, Theorems \ref{Main_theorem_about_approximation_of_IS_method}, \ref{Second_theorem_about_approximation_properties_of_SI_method}, deals with the SI-method for equation \eqref{first_eq} subjected to an initial condition, and provides a priori error estimates for the case.

The paper is organized as follows. In the beginning of Section \ref{Section_label_SI_method_for_ivp} we introduce the SI-method for solving initial value problems associated
with equation \eqref{first_eq}; the rest part of the section is devoted to a thorough investigation of the method's approximation properties, which are formulated as Theorems \ref{Main_theorem_about_approximation_of_IS_method} and \ref{Second_theorem_about_approximation_properties_of_SI_method}. The SI-method for solving boundary value problems \eqref{first_eq} \eqref{boundary_cond} is the main focus of Section \ref{Section_label_SI_method_for_bvp}, where we describe a {\it single } and {\it multiple shooting } versions of the method. We apply the SI-method to the Troesch's equation subjected to both initial and boundary conditions and discuss the results in Section \ref{Section_label_numerical_examples}. Section \ref{Section_label_conclusions} contains conclusions.

\section{Straight-Inverse method for solving IVPs for the second order differential equations}\label{Section_label_SI_method_for_ivp}

\subsection{Step functions}\label{Section_step_functions}
Before proceeding any further with the description of the SI-method, we need to introduce a pair of, so called, \textit{step functions} $U(s)$ and $V(s),$ which play an important role in the method's framework. It is worth mentioning that, in principle, the functions can be chosen in a multiple different ways, resulting in different implementations of the method. For the sake of simplicity, below we give a very concrete definition of the step functions and stick to it throughout the rest of the paper.

\begin{definition}
We define step function $U(s)$ to be the solution for the initial value problem
\begin{equation}\label{gen_problem_straight}
\frac{d^{2}U(s)}{d s^{2}} = (As + B)U(s),  \quad U(0) = D;\quad U^{\prime}(0) = C, \quad A, B, C, D, s \in \mathbb{R}.
\end{equation}
\end{definition}

\begin{definition}
We define step function $V(s)$ to be the solution for the initial value problem
\begin{equation}\label{gen_problem_inverse}
\frac{d^{2} V(s)}{d s^{2}} = (\bar{A} s + \bar{B})\frac{d V(s)}{d s},\quad V(0) = \bar{D},\quad V^{\prime}(0) = \bar{C}, \quad \bar{A}, \bar{B}, \bar{C}, \bar{D}, s\in \mathbb{R},
\end{equation}
\end{definition}

One might notice that if
\begin{equation}
    \begin{split}
        A = & N^{\prime}_{u}(u_{a}, a)u^{\prime}_{a} + N^{\prime}_{x}(u_{a}, a),\\
        B = & N(u_{a}, a),\\
        C = & u^{\prime}_{a},\\
        D = & u_{a},
    \end{split}
\end{equation}
then the function $U(x - a)$ coincides with the solution to the linearization of equation \eqref{first_eq}, supplemented with the initial conditions
\begin{equation}\label{initial_conditions}
    u(a) = u_{a} \in \mathbb{R},\quad u^{\prime}(a) = u^{\prime}_{a}\in \mathbb{R}.
\end{equation}

The meaning of step function $V(s)$ becomes more clear in the light of the statement below.
\begin{lemma}\label{Lemma_x_is_a_solution_to_inverse_problem}
  Let $u(x)$ be the unique solution to IVP \eqref{first_eq}, \eqref{initial_conditions} on $[a, b],$ $a < b$. If $u^{\prime}_{a} > 0$ $(u^{\prime}_{a} < 0)$ and $u(x)$ is monotone on $[a, b]$ then the inverse function $x(\cdot) = u^{-1}(\cdot)$ is the unique solution to the IVP
\begin{eqnarray}
\frac{d^{2} x(u)}{d u^{2}} = - N\left(u, x(u)\right)u \left(\frac{d x(u)}{d u}\right)^{3},  \label{first_eq_inverse} \\
x(u_{a}) = a, \quad x^{\prime}(u_{a}) = x^{\prime}_{a} = 1/u^{\prime}_{a}, \label{boundary_cond_inverse}
\end{eqnarray}
on $[u_{a}, u(b)]$ $([u(b), u_{a}])$.
\end{lemma}

Now if we assume that
\begin{equation}
\begin{split}
\bar{A} = & - \left((N^{\prime}_{u}(u_{a}, a) + N_{x}^{\prime}(u_{a}, a)x^{\prime}_{a})u_{a} +  N(u_{a}, a)\right)\left(x^{\prime}_{a}\right)^{2} + 2 \left(N(u_{a}, a)u_{a}\right)^{2}\left(x^{\prime}_{a}\right)^{4}, \\
\bar{B} = & - N(u_{a}, a)u_{a}\left(x^{\prime}_{a}\right)^{2},\\
\bar{C} = & x^{\prime}_{a},\\
\bar{D} = & a,\\
\end{split}
\end{equation}
then the function $V(u - u_{a})$ is nothing else but the solution to the linearization of equation \eqref{first_eq_inverse} subjected to initial conditions \eqref{boundary_cond_inverse}.

\subsection{Description of the SI-method for solving IVPs}

Let $\Omega(h)$ denotes an ordered set of quadruples of the form
\begin{equation}\label{Definitions_of_omega_mesh}
  \Omega(h) = \{\left(u^{\prime}_{i}, x^{\prime}_{i}, u_{i}, x_{i}\right), \; u^{\prime}_{i} \stackrel{def}{=} 1/x^{\prime}_{i}, i=0,1,2,\ldots\}
\end{equation}
with elements defined by means of the following chain of recurrence equalities:

\begin{equation} \label{SI_def_first}
x_{0} = a,\quad u_{0} = u_{a}, \quad u^{\prime}_{0} = u_{a}^{\prime};
\end{equation}
if $|u^{\prime}_{i-1}| \leq 1$ then
\begin{equation}\label{SI_def_second}
\begin{split}
x_{i} = & x_{i-1} + h,\\
u_{i} =  & U(A_{i-1}, B_{i-1}, C_{i-1}, D_{i-1}, h), \\
u^{\prime}_{i} = & U^{\prime}_{h}(A_{i-1}, B_{i-1}, C_{i-1}, D_{i-1}, h),\\
A_{i} = & N^{\prime}_{u}(u_{i}, x_{i})u^{\prime}_{i} + N_{x}^{\prime}(u_{i}, x_{i}),\\
B_{i} = & N(u_{i}, x_{i}),\\
C_{i} = & u^{\prime}_{i},\\
D_{i} = & u_{i},\\
\end{split}
\end{equation}
otherwise, if $|u^{\prime}_{i-1}| > 1$
\begin{equation}\label{SI_def_third}
\begin{split}
x_{i} = & V(\bar{A}_{i-1}, \bar{B}_{i-1}, \bar{C}_{i-1}, \bar{D}_{i-1}, h_{i}^{\ast}),\\
u_{i} = & u_{i-1} + h_{i}^{\ast},\\
x^{\prime}_{i} = & V^{\prime}_{h}(\bar{A}_{i-1}, \bar{B}_{i-1}, \bar{C}_{i-1}, \bar{D}_{i-1}, h_{i}^{\ast}),\\
\bar{A}_{i} = & - \left((N^{\prime}_{u}(u_{i}, x_{i}) + N_{x}^{\prime}(u_{i}, x_{i})x_{i}^{\prime})u_{i} \right. + \\ + & N(u_{i}, x_{i})\big)\left(x_{i}^{\prime}\right)^{2} + 2 \left(N(u_{i}, x_{i})u_{i}\right)^{2}\left(x_{i}^{\prime}\right)^{4}, \\
\bar{B}_{i} = & - N(u_{i}, x_{i})u_{i}\left(x_{i}^{\prime}\right)^{2},\\
\bar{C}_{i} = & x_{i}^{\prime},\\
\bar{D}_{i} = & x_{i},\\
h_{i}^{\ast} = & \sign(x^{\prime}_{i-1}) h,
\end{split}
\end{equation}
where
 $h$ --- some fixed positive real number hereinafter referenced to as a \textit{step size} of the SI-method. Formulas \eqref{SI_def_second} can be interpreted as a "straight" phase of the method, since they deal with the "straight" problem \eqref{first_eq}, \eqref{initial_conditions}, whereas formulas \eqref{SI_def_third} describe the method's "inverse" phase, dealing with the "inverse" problem
 \eqref{first_eq_inverse}, \eqref{boundary_cond_inverse}.

Ordered set $\Omega(h)$ \eqref{Definitions_of_omega_mesh} will be referenced to as a \textit{mesh} of the SI-method that corresponds to IVP \eqref{first_eq}, \eqref{initial_conditions}. From the recurrence formulas \eqref{SI_def_first}, \eqref{SI_def_second},
\eqref{SI_def_third} it follows that if function $N(u, x)$ belongs to $C^{1}(\mathbb{R}\times[0, +\infty))$ then the mesh $\Omega(h)$ contains
infinite number of elements, i.e. the recurrence process of calculating quadruples $\left(u^{\prime}_{i}, x^{\prime}_{i}, u_{i}, x_{i}\right)$ can be continued infinitely long.
In the light of this, a reasonable question arises: whether the mesh $\Omega(h)$ (which is infinite) have something to do with the exact solution $u(x)$
of the Cauchy problem \eqref{first_eq}, \eqref{initial_conditions} (which might exist only on some finite subinterval of $[a, +\infty)$) and, if yes,
what approximation properties does the mesh possess with respect to the exact solution? The question is addressed in the paragraph below.

\subsection{Error analysis.}

Although the SI-method, introduced above, is applicable to a class of IVPs associated with equation \eqref{first_eq}, for the sake of simplicity,
the main theoretical results, revealing approximation properties of the method, are stated and proved for a more narrow set of problems, as it can be seen from the theorems below.
\begin{theorem}\label{Main_theorem_about_approximation_of_IS_method}
Let the nonlinear function $N(u, x)$ be independent on $x$, i.e.
\begin{equation}\label{equation_n_independent_on_x}
  \frac{\partial N(u, x)}{\partial x} \equiv 0,
\end{equation}
and
\begin{equation}\label{equation_n_continuously_differentiable}
  N(u) \equiv N(u, x) \in C^{2}([u_{a}, +\infty)),
\end{equation}
\begin{equation}\label{equation_n_positive}
  N(u) > 0,\; N^{\prime}(u) \geq 0 \; \forall u\in (u_{a}, +\infty),
\end{equation}
\begin{equation}\label{restrictions_on_initial_conditions}
  u_{a} \geq 0,\; 0 < u^{\prime}_{a} \leq 1,
\end{equation}
\begin{equation}\label{equation_n_integral_increases}
  \lim\limits_{u \rightarrow +\infty} \frac{1}{u^{2 + \lambda}}\int\limits_{u_{a}}^{u}N(\xi)\xi d\xi > 0,
\end{equation}
for some $\lambda > 0.$

If
\begin{equation}\label{theorem_conditions_first_main_restriction_on_h}
  0 < h < \min\left\{1, \sqrt{\frac{\varepsilon}{P}}, \frac{\varepsilon}{ L_{0}M_{0}}\right\}
\end{equation}
then there exists an integer $i^{\ast} > 0,$ such that
\begin{equation}\label{thoerem_conditions_i_ast_definition}
  u^{\prime}_{i} \leq 1, \forall i \in \overline{0, i^{\ast} - 1},\; u^{\prime}_{i^{\ast}} > 1
\end{equation}
and the estimate holds true
\begin{equation}\label{theorem_conditions_first_estimation}
  \max\{|u(x_{i}) - u_{i}|, |u^{\prime}(x_{i}) - u^{\prime}_{i}|\} \leq h^{2}P,\; \forall i\in \overline{0, i^{\ast}},
\end{equation}
where $x_{i},$ $u_{i},$ $u_{i}^{\prime}$ are calculated according to formulas \eqref{SI_def_first}, \eqref{SI_def_second}; $u(x)$ is the solution of IVP \eqref{first_eq}, \eqref{initial_conditions};
\begin{equation}\label{theorem_condition_p_ast_definition}
  P = M_{0}M_{1}\left( L_{2}M_{1} +  L_{1} L_{0}  M_{0}\right)\times
\end{equation}
$$\times\frac{\exp\left((S^{\ast}+1)\max\{1,  L_{0} +  L_{1}\hat{M}_{1}\} + S^{\ast}M_{0}(2 L_{1} +  L_{2}\hat{M}_{1})\right)}
  {2\left( L_{1}M_{0} + \max\{1,  L_{0}\}\right)},$$
\begin{equation}\label{Proof_L_constant_definition}
   L_{i} = \max\limits_{|u| < M_{0} + \varepsilon} \left|N^{(i)}(u)\right|,\; i= 0,1,2,
\end{equation}
\begin{equation}\label{Theorem_conditions_M_constant_definition}
  M_{0} = \frac{1}{2}S^{\ast}(M_{1} - u_{a}^{\prime}), \; M_{1} = 1+3\varepsilon,\; \hat{M}_{1} = M_{1} + \varepsilon,
\end{equation}
\begin{equation}\label{theorem_conditions_S_ast_definition}
  S^{\ast} = \lim\limits_{u \rightarrow +\infty}\int\limits_{u_{a}}^{u}\frac{d\eta}{\sqrt{(u^{\prime}_{a})^{2} + 2\int\limits_{u_{a}}^{\eta}N(\xi)\xi d\xi}}
  \footnote{The existance of the limit follows from condition \eqref{equation_n_integral_increases},}
\end{equation}
and $\varepsilon$ denotes an arbitrary positive parameter.
 In addition to that, an auxiliary estimate holds true
\begin{equation}\label{Theorem_aux_estimete_magnitude}
  u(x) \leq M_{0}, \forall x\in [x_{0}, x_{i^{\ast}} + h].
\end{equation}
\end{theorem}

\begin{proof}
   As it is stated in Lemma \ref{Lemma_x_is_a_solution_to_inverse_problem}, the function $x(u),$
   which is (by definition) inverse of the exact solution $u(x),$ should be the solution to IVP \eqref{first_eq_inverse}, \eqref{boundary_cond_inverse}.
   Under the assumptions of the theorem, equation \eqref{first_eq_inverse} becomes a partial case of the well known Bernoulli equation, which allows us to express
   the solution $x(u)$ in the closed form (see, for example, \cite{zaitsev2002handbook}):
   \begin{equation}\label{inverse_solution_closed_form}
     x(u) = a + \int\limits_{u_{a}}^{u}\frac{d\eta}{\sqrt{(u^{\prime}_{a})^{2} + 2\int\limits_{u_{a}}^{\eta}N(\xi)\xi d\xi}}.
   \end{equation}
    From \eqref{equation_n_independent_on_x}, \eqref{equation_n_continuously_differentiable}, \eqref{equation_n_positive}, \eqref{restrictions_on_initial_conditions} and the Picard-Lindel\"{o}f theorem
   (see \cite[p. 38]{Teschl_ODE_and_DS}) it follows that function $x(u)$ \eqref{inverse_solution_closed_form} belongs to $C^{3}([u_{a}, +\infty))$ and is the unique solution to the IVP
   \eqref{first_eq_inverse}, \eqref{boundary_cond_inverse} on $[u_{a}, +\infty).$

   Using inequalities \eqref{equation_n_positive}, \eqref{restrictions_on_initial_conditions}, \eqref{equation_n_integral_increases} and the
   \textit{Limit Comparison Theorem for Improper Integrals}, from \eqref{inverse_solution_closed_form} we can easily derive
   that $x(u)$ is a monotonically increasing function on $[u_{a}, +\infty)$ with bounded range:
   $$[u_{a}, +\infty) \stackrel{x(\cdot)}{\Rightarrow} [a, S), \; S = \lim\limits_{u \rightarrow +\infty} x(u) < +\infty. $$
   The letter fact means that its inverse, $u(x),$ exists on $\bigl[a, S\bigr)$ and is the
   unique solution to IVP \eqref{first_eq}, \eqref{initial_conditions} on the segment. Furthermore, function $u(x)$ is monotonically increasing on $[a, S)$ and, taking into account condition $u_{a} \geq 0$ \eqref{restrictions_on_initial_conditions}, positive on $(a, S),$ i.e.
   \begin{equation}\label{monotonicity_of_u}
     0 < u(\xi_{1}) < u(\xi_{2}),\; \forall \xi_{1}, \xi_{2} \in (a, S) : \xi_{1} < \xi_{2}.
   \end{equation}

   As it follows from representation
   \eqref{inverse_solution_closed_form}, conditions \eqref{equation_n_positive}, \eqref{restrictions_on_initial_conditions}, \eqref{equation_n_integral_increases} also mean that $x^{\prime}(u)$ is
   a positive, monotonically decreasing function on $[u_{a}, +\infty),$ which tends to $0$ as $u$ tends to $+\infty.$ Consequently, $u^{\prime}(x)$ is a positive, monotonically increasing function on $[a, S)$, which tends to $+\infty$ as $x$ tends to $S:$
   \begin{equation}\label{deriv_of_straight_tends_to_infinity}
   0 < u^{\prime}(\xi_{1}) < u^{\prime}(\xi_{2}), \; \forall \xi_{1}, \xi_{2}\in [a, S) : \xi_{1} < \xi_{2},\; \lim_{x \rightarrow +S} u^{\prime}(x) = +\infty.
   \end{equation}

   From \eqref{deriv_of_straight_tends_to_infinity} it follows that for each $\delta \geq u_{a}^{\prime},$ there exists a unique $x_{\delta} \in [a, S)$ such that
   $$u^{\prime}(x_{\delta}) = \delta.$$ The latter, in conjunction with the fact that function $u^{\prime}(x)$ is convex on $[a, S),$\footnote{I.e. $u^{\prime\prime\prime}(x) = N^{\prime}(u(x))u(x) + N(u(x))u^{\prime}(x) > 0, \; \forall x\in (a, S).$}
   allows us to establish the inequality (see Fig. \ref{fig:M1})
   \begin{equation}\label{inequality_for_the_norm_of_u_through_value_of_u_prime}
     \max\limits_{x\in [a, S) : u^{\prime}(x) \leq \delta}|u(x)| = \int\limits_{a}^{x_{\delta}}u^{\prime}(\xi)d\xi\leq
     \frac{1}{2}(x_{\delta} - a)(\delta - u_{a}^{\prime}) <
     \frac{1}{2}(S - a)(\delta - u_{a}^{\prime}),
   \end{equation}
   which is of crucial importance for the rest of the proof.

    \begin{figure}
    \centering
    \begin{tikzpicture}
    \draw[<->] (6,0) node[below]{$x$} -- (0,0) --
    (0,6) node[left]{$u$};
    \draw (5.05, -0.1) node[below]{$S$} -- (5.05, 6);
    \draw (0, 0.5) node[left]{$A(a,u^{\prime}_{a})$} -- (6, 0.5);
    \draw (6.3, 3.75) node[left]{$B(x_{\delta},\delta)$};
    \draw (6.4, 0.75) node[left]{$C(x_{\delta},u^{\prime}_{a})$};
    \draw (0, 3.5) node[left]{$\delta$} -- (6, 3.5);
    \draw (4.6, -0.1) node[below]{$x_{\delta}$} -- (4.6, 6);
    \draw [dashed] (0, 0.5) -- (4.6, 3.5);
    \draw[very thick] (0,0.5) to [out=0,in=-110] (4.6, 3.5) to [out = 70, in = -90] (5,5.5);
    \path [fill=gray] (4.6, 3.5) -- (4.6,0.5) -- (0,0.5) to [out=0,in=-110] (4.6, 3.5);
    \end{tikzpicture}
    \caption{Graph of $u=u^{\prime}(x)$ (solid line). Area of the shaded region is equal to $\int\limits_{a}^{x_{\delta}}u^{\prime}(\xi)d\xi$, which is, apparently, less or equal to the area of
    $\triangle ABC,$ which, in turn, is equal to $\frac{1}{2}(x_{\delta} - a)(\delta - u_{a}^{\prime}).$ } \label{fig:M1}
    \end{figure}
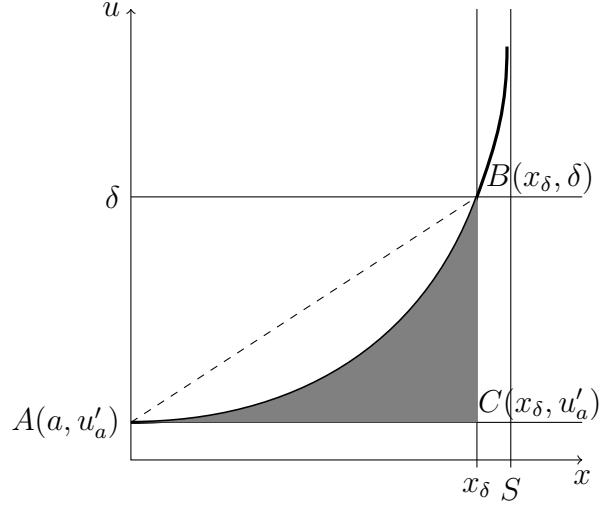

   Using notation $$ e_{i} = \max\{|u(x_{i}) - u_{i}|, |u^{\prime}(x_{i}) - u^{\prime}_{i}|\}, $$
   we can estimate $e_{i}$ from the system of differential equation
   \begin{equation}\label{Proof_z_i_general_equation}
     \dot{Z_{i}}(x) =
    \left[\begin{array}{cc}
    0 & 1\\
    N(u_{i-1}) + N^{\prime}(u_{i-1})u^{\prime}_{i-1}(x- x_{i-1}) & 0
    \end{array}\right]Z_{i}(x) +
    \left[\begin{array}{c}
    0  \\
    F_{i}(x)u(x)
    \end{array}\right],
   \end{equation}
   $$ x\in [x_{i-1}, x_{i}],,\; Z_{i}(x_{i}) = Z_{i-1}(x_{i}),$$
    where
    $$F_{i}(x) = N(u(x)) - N(u_{i-1}) - N^{\prime}(u_{i-1})u^{\prime}_{i-1}(x- x_{i-1}),$$
    $$Z_{i}(x) = \left[
                \begin{array}{c}
                  z_{i}(x) \\
                  z^{\prime}_{i}(x) \\
                \end{array}
              \right], \;\; z_{i}(x) = u(x) - U(N^{\prime}(u_{i-1})u^{\prime}_{i-1}, N(u_{i-1}), u^{\prime}_{i-1}, u_{i-1}, x - x_{i-1}), \; i=1,2,\ldots,
    $$
    \begin{equation}\label{Proof_Z_zero_definition}
      Z_{0}(x)\equiv 0.
    \end{equation}
Indeed,
$$e_{i} \leq \|Z_{i}(x)\|_{[x_{i-1}, x_{i}]} \stackrel{def}{=} \max\limits_{x\in [x_{i-1}, x_{i}]}\|Z_{i}(x)\| = \max\limits_{x\in [x_{i-1}, x_{i}]}\max\{z_{i}(x), z^{\prime}(x)\}.$$
We are not going to estimate $e_{i}$ for all integer $i,$ but only for those satisfying inequality $i \leq i^{\ast},$ where $i^{\ast}$
is defined in \eqref{thoerem_conditions_i_ast_definition}. However, at this point, the very existence of such an integer value $i^{\ast}$ is yet to be proved.

To prove that $i^{\ast}$ exists, let us fix some arbitrary $\varepsilon > 0$ and make an assumption that
\begin{equation}\label{Proof_assumption_about_deviation}
  \|Z_{j}(x)\|_{[x_{j-1}, x_{j}]} < \varepsilon,\; j = 1,2,\ldots, j_{3\varepsilon};\; j_{3\varepsilon} = \max \{i \;|\; u^{\prime}(x_{i}) \leq M_{1} = 1+3\varepsilon\},
\end{equation}
for $h$ sufficiently small. For the given $\varepsilon,$ we also consider constants $ L_{i},$ $i=0,1,2$ defined in \eqref{Proof_L_constant_definition}, \eqref{Theorem_conditions_M_constant_definition}, keeping in mind that, according to inequality \eqref{inequality_for_the_norm_of_u_through_value_of_u_prime},
\begin{equation}\label{Proof_M_constant_definition}
  M_{0} = \frac{1}{2}(S - a)(M_{1} - u_{a}^{\prime}) > \max\limits_{x\in [a; S) : u^{\prime}(x) \leq M_{1}} |u(x)|.
\end{equation}

Now, requiring that
\begin{equation}\label{Proof_second_upper_bound_restriction_for_h}
  h < \frac{\varepsilon}{ L_{0} M_{0}},
\end{equation}
we can easily prove that $i^{\ast}$ exists and $x_{i^{\ast}}$ belongs to $(a, u^{\prime -1}(1+2\varepsilon)).$ Indeed, if there exists $x_{j}\in (a, u^{\prime -1}(1 + 2\varepsilon)],$
such that $u^{\prime}_{j} > 1,$ then we can put $$i^{\ast} = \min\left\{j \in \overline{1, j_{3\varepsilon}} \; | \; u^{\prime}_{j}> 1\right\}.$$
If this is not the case at least for a single $h$ satisfying \eqref{Proof_second_upper_bound_restriction_for_h}, then, from \eqref{Proof_second_upper_bound_restriction_for_h} it follows that
$$h < u^{\prime -1}(1+2\varepsilon) - u^{\prime -1}(1 + \varepsilon) \geq \frac{\varepsilon}{\max\limits_{x\in (a, u^{\prime -1}(1+2\varepsilon))}|u^{\prime\prime}(x)|}\geq \frac{\varepsilon}{L_{0}M_{0}},$$
which, in turn, means that there exists at least one $x_{j}$ belonging to the interval
$$(u^{\prime -1}(1 + \varepsilon), u^{\prime -1}(1 + 2\varepsilon)].$$ Taking into account assumption \eqref{Proof_assumption_about_deviation},
the latter fact yields us
$$u^{\prime}_{j} = u^{\prime}(x_{j}) + u^{\prime}_{j} - u^{\prime}(x_{j}) > 1 + \varepsilon - |u^{\prime}_{j} - u^{\prime}(x_{j})| > 1 + \varepsilon - \varepsilon = 1$$
and, consequently, we get a contradiction.

By a similar reasoning, we can easily prove that $x_{i^{\ast}} + h < u^{\prime -1}(1+3\varepsilon)$, which, being combined with \eqref{Proof_M_constant_definition}, gives us auxiliary estimate \eqref{Theorem_aux_estimete_magnitude}.

Using constants  \eqref{Proof_L_constant_definition}, \eqref{Theorem_conditions_M_constant_definition} and assumption \eqref{Proof_assumption_about_deviation}, from \eqref{Proof_z_i_general_equation} we can derive recurrent estimates for $\|Z_{i}(x)\|_{[x_{i-1}, x_{i}]},$ $i = 1,2,\ldots, j_{3\varepsilon}$
in the following way:
\begin{equation}\label{Proof_general_inequality_for_Z_i_norm}
  \|Z_{i}(x)\| \leq \left(1+h Q\right)\|Z_{i-1}(x)\|_{[x_{i-2}, x_{i-1}]} + E\int\limits_{x_{i-1}}^{x}\|Z_{i}(\xi)\|d\xi + h^{3} K,\; \forall x\in [x_{i-1}, x_{i}],\;
\end{equation}
$$Q = M_{0}\left( L_{1} + h( L_{1} +  L_{2}\hat{M}_{1})\right),\; E = \max\left\{1,  L_{0} +  L_{1}\hat{M}_{1}h\right\},$$
$$ K = \frac{M_{0}M_{1}}{2} \left( L_{2}M_{1} +  L_{1} L_{0} M_{0}\right), \; \hat{M}_{1} = M_{1} + \varepsilon \geq u_{i}^{\prime}\footnote{The inequality follows from the definition of $M_{1}$ \eqref{Theorem_conditions_M_constant_definition} and assumption \eqref{Proof_assumption_about_deviation}.}$$
$$i = 1,2,\ldots, j_{3\varepsilon},\; \|Z_{0}(x)\|_{[x_{-1}, x_{0}]} \equiv 0.$$

Applying the Gronwall's inequality (see, for example, \cite[42]{Teschl_ODE_and_DS}) to \eqref{Proof_general_inequality_for_Z_i_norm} we get
\begin{equation}\label{Proof_general_inequality_for_Z_i_with_Gronwall}
   \|Z_{i}(x)\|_{[x_{i-1},x_{i}]} \leq \left(\left(1+h Q\right)\|Z_{i-1}(x)\|_{[x_{i-2}, x_{i-1}]} + h^{3} K\right)\exp\left(h E\right), \; i=1,2,\ldots, j_{3\varepsilon}.
\end{equation}

Inequality \eqref{Proof_general_inequality_for_Z_i_with_Gronwall}, in conjunction with \eqref{Proof_Z_zero_definition}, yields us the estimate
\begin{equation}\label{Proof_Z_i_estimation}
  \|Z_{i}(x)\|_{[x_{i-1},x_{i}]} \leq h^{3} K \sum\limits_{j=1}^{j=i}(1+hQ)^{j-1}\exp\left(j h E\right) =
\end{equation}
$$ = h^{3} K\exp\left(h E\right)\frac{(1+h Q)^{i}\exp\left(ihE\right) - 1}{(1+h Q)\exp\left(h E\right) - 1} \leq h^{3} K\exp\left(h E\right)\frac{\exp\left(S^{\ast}\left(E+ Q\right)\right) - 1}{(1+h Q)\exp\left(h E\right) - 1} \leq $$
$$ \leq h^{2} K\frac{\exp\left((S^{\ast}+1)\overline{E} + S^{\ast}M_{0}(2 L_{1} +  L_{2}\hat{M}_{1})\right)}{ L_{1}M_{0} + \underline{E}} = h^{2} P,\; i = 1, 2, \ldots, j_{3\varepsilon},$$
where
$$ \overline{E} = \max\{1,  L_{0} +  L_{1}\hat{M}_{1}\},\; \underline{E} = \max\{1,  L_{0}\},\; S^{\ast} = S -a.$$

The last inequality in \eqref{Proof_Z_i_estimation} holds true under the assumption
\begin{equation}\label{Proof_assumption_h_less_or_equal_than_1}
  0 \leq h \leq 1,
\end{equation}
which we accept from now on.

With estimate \eqref{Proof_Z_i_estimation} in hands, we can focus on proving inequalities \eqref{Proof_assumption_about_deviation}, which have been embraced as an assumption until now. Going back to inequality \eqref{Proof_general_inequality_for_Z_i_norm}, it is important to mention that to derive it for each particular $i=1,2,\ldots, j_{3\varepsilon},$ we need to use assumption \eqref{Proof_assumption_about_deviation} for $j = i-1$ only. Besides that, inequality \eqref{Proof_general_inequality_for_Z_i_norm} (and, consequently, inequality \eqref{Proof_Z_i_estimation}) for $i = 1$ does not rely upon \eqref{Proof_assumption_about_deviation} at all. From \eqref{Proof_Z_i_estimation} it follows that if we require
\begin{equation}\label{Proof_first_upper_bound_restriction_for_h}
  h < \sqrt{\frac{\varepsilon}{P}},
\end{equation}
then inequality \eqref{Proof_assumption_about_deviation} for $j=1$ holds true (i.e., immediately follows from \eqref{Proof_Z_i_estimation}). The latter automatically implies inequality \eqref{Proof_Z_i_estimation} for $i=2,$ which, together with \eqref{Proof_first_upper_bound_restriction_for_h}, yields the fulfillment of \eqref{Proof_assumption_about_deviation} for $j=2\ldots$ Apparently, using the method of mathematical induction, we can easily proof that
under condition \eqref{theorem_conditions_first_main_restriction_on_h} all the inequalities \eqref{Proof_assumption_about_deviation} hold true.

This concludes proof of the theorem, stating existence of the index $i^{\ast},$ satisfying conditions \eqref{thoerem_conditions_i_ast_definition}, and fulfillment of estimates \eqref{theorem_conditions_first_estimation}, \eqref{Theorem_aux_estimete_magnitude}, provided that $h$ satisfies \eqref{theorem_conditions_first_main_restriction_on_h}.
\end{proof}

The estimates given in Theorem \ref{Main_theorem_about_approximation_of_IS_method} are quite rough. In part, this is because of the roughness of estimate \eqref{Proof_M_constant_definition}.
The latter can be improved, as it is proposed in the remarks below.

\begin{remark}\label{Remark_to_main_theorem}
The estimates of Theorem \ref{Main_theorem_about_approximation_of_IS_method} remain valid and can be improved if the constant $M_{0}$ defined by
 formula \eqref{Theorem_conditions_M_constant_definition} is substituted by the one defined as
\begin{equation}\label{Remark_1_re-definition_of_m_constant}
  M_{0} = \Phi^{-1}\left(\frac{1}{2}\left(M_{1}^{2} - \left(u^{\prime}_{a}\right)^{2}\right)\right), \; \Phi(u) = \int\limits_{u_{a}}^{u}N(\xi)\xi d\xi,\; M_{1} = 1+3\varepsilon,
\end{equation}
and the constant $P$ defined at \eqref{theorem_condition_p_ast_definition} is treated as a function of $h$ defined as
\begin{equation}\label{Remark_1_re-definition_of_P_constant}
  P(h) = hK\exp\left(h E\right)\frac{\exp\left(S^{\ast}\left(E+ Q\right)\right) - 1}{\left(1+h Q\right)\exp\left(h E\right) - 1},
\end{equation}
where
$$Q = M_{0}\left( L_{1} + h( L_{1} +  L_{2}\hat{M}_{1})\right),\; E = \max\left\{1,  L_{0} +  L_{1}\hat{M}_{1}h\right\},\; K = \frac{M_{0}M_{1}}{2} \left( L_{2}M_{1} +  L_{1} L_{0} M_{0}\right).$$
\end{remark}

\begin{remark}\label{Remark_about_estimate_for_S}
  If the IVP \eqref{first_eq}, \eqref{initial_conditions} is considered on some finite interval, i.e. $b < +\infty$ then the constant $S^{\ast}$ defined in
  \eqref{theorem_conditions_S_ast_definition} can be substituted with
  \begin{equation}\label{Remark_2_re_definition_om_S}
    S^{\ast} = \min\left\{b, \lim\limits_{u \rightarrow +\infty}\int\limits_{u_{a}}^{u}\left((u^{\prime}_{a})^{2} + 2\int\limits_{u_{a}}^{\eta}N(\xi)\xi d\xi\right)^{-\frac{1}{2}}d\eta\right\}
  \end{equation}
  in order to make error estimates of Theorem \ref{Main_theorem_about_approximation_of_IS_method} more precise.
\end{remark}

As one might notice, Theorem \ref{Main_theorem_about_approximation_of_IS_method} is concerned with the "straight" phase of the SI-method, which is described by formulas \eqref{SI_def_first}, \eqref{SI_def_second}. The "inverse" phase of the method is the main focus of the theorem below.

\begin{theorem}{\label{Second_theorem_about_approximation_properties_of_SI_method}}
Let the conditions of Theorem \ref{Main_theorem_about_approximation_of_IS_method} are fulfilled and the notations $i^{\ast}, x_{i^{\ast}}, u_{i^{\ast}}, u_{i^{\ast}}^{\prime},$ $M_{0},$ $P,$ $L_{0,1,2},$ $\varepsilon$ keep their meaning (from the Theorem \ref{Main_theorem_about_approximation_of_IS_method}).
If, additionally,
\begin{equation}\label{Second_theorem_definition_of_N_bold}
  \mathcal{N}(u) \stackrel{def}{=}\left\{
                                      \begin{array}{cc}
                                        N(u)u & u \geq u_{i^{\ast}}, \\
                                        0 & u < u_{i^{\ast}}, \\
                                      \end{array}
                                    \right.
\end{equation}
\begin{equation}\label{theorem_conditions_definition_of_mu}
  \mu = \sup\limits_{u\in [u_{i^{\ast}}, +\infty)} \frac{\mathcal{N}(u)}{1+ \int\limits_{u_{i^{\ast}}}^{u}\mathcal{N}(\xi) d\xi} < +\infty,
\end{equation}
\begin{equation}\label{Theorem_1_epsilon_definition}
 0 < \varepsilon < \frac{1}{2},
\end{equation}
\begin{equation}\label{theorem_conditions_second_main_restriction_on_h}
  h < \min\left\{\frac{1-2\varepsilon}{P}, \frac{1}{3\mu}\right\},
\end{equation}
\begin{equation}\label{equation_n_derivative_positive}
  \mathcal{N}(u) \in C^{3}([u_{i^{\ast}}, +\infty]),\; \mathcal{N}^{(k)}(u) \geq 0, \; k=0,1,2,3, \; \forall u\in [u_{i^{\ast}}, +\infty],
\end{equation}
then the following estimates hold true for all $i \geq i^{\ast}$
\begin{equation}\label{theorem_two_x_i_prime_relains_below_1}
  x_{i}^{\prime} < 1,
\end{equation}
\begin{equation}\label{theorem_statement_about_x_approximation}
  |x_{i} - x(u_{i})| \leq \frac{Ph^{2}}{\tau(h)}+
\end{equation}
$$+ h^{2}\int\limits_{u_{i^{\ast}}}^{u_{i}}\exp\left(\int\limits_{u_{i^{\ast}}}^{\eta}\Lambda(\zeta, h)d\zeta\right)\!\!\!\left(\left(\frac{L_{0}M_{0}}{\tau(h)^{3}} + \frac{1}{1-Ph^{2}}\right)P + \frac{1}{2}\int\limits_{u_{i^{\ast}}}^{\eta}\mathcal{T}(\zeta, h)d\zeta \right)d\eta,$$
\begin{equation}\label{theorem_statement_about_x_prime_approximation}
  |x^{\prime}_{i} - x^{\prime}(u_{i})| \leq h^{2}\exp\left(\int\limits_{u_{i^{\ast}}}^{u_{i}}\Lambda(\zeta, h)d\zeta\right)\!\!\!\left(\left(\frac{L_{0}M_{0}}{\tau(h)^{3}} + \frac{1}{1-Ph^{2}}\right)P + \frac{1}{2}\int\limits_{u_{i^{\ast}}}^{u_{i}}\mathcal{T}(\zeta, h)d\zeta \right),
\end{equation}
where $x_{i}, u_{i}, x_{i}^{\prime}$ for $i > i^{\ast}$ are calculated according to formulas \eqref{SI_def_third}, $x(\cdot) \overset{def}{=} u^{-1}(\cdot),$
$$\tau(h) = 1 -  L_{0}M_{0} h - Ph^{2},$$
\begin{equation}\label{Exponenta_term_expression}
  \Lambda(\zeta, h) = 2\left(\frac{\mathcal{N}(\zeta) + h\mathcal{N}^{\prime}(\zeta)}{\left(\Upsilon_{2}(\zeta - 2h)\right)^{\frac{1}{2}}} + \frac{2h\left(\mathcal{N}(\zeta)\right)^{2}}{\left(\Upsilon_{2}(\zeta -2h)\right)^{\frac{3}{2}}}\right)\frac{d\zeta}{\sqrt{\Upsilon_{1}(\zeta, h)}},
\end{equation}
\begin{equation}\label{Second_derivative_term_expression}
  \mathcal{T}(\zeta, h) = \left(\frac{\mathcal{N}^{\prime\prime}(\zeta)}{\Upsilon_{1}(\zeta - h, h)} + \frac{6\mathcal{N}^{\prime}(\zeta)\mathcal{N}(\zeta)}{\left(\Upsilon_{1}(\zeta-h,h)\right)^{2}}
  +\frac{8\left(\mathcal{N}(\zeta)\right)^{3}}{\left(\Upsilon_{1}(\zeta-h, h)\right)^{3}}\right)\frac{1}{\sqrt{\Upsilon_{1}(\zeta,h)}},
\end{equation}
\begin{equation}\label{Second_theorem_Upsilon_1_2_definition}
  \Upsilon_{1}(u, h) = \tau\left(h\right)^{2} + 2\int\limits_{u_{i^{\ast}}}^{u}\mathcal{N}(\xi)d\xi, \; \Upsilon_{2}(u) = 1 + \int\limits_{u_{i^{\ast}}}^{u}\mathcal{N}(\xi)d\xi.
\end{equation}
 \end{theorem}

To prove Theorem \ref{Second_theorem_about_approximation_properties_of_SI_method} we will need an auxiliary statement below.
\begin{lemma}\label{Proof_lemma_about_stability_of_V_function}
Let function $\mathcal{N}(u),$ defined as \eqref{Second_theorem_definition_of_N_bold} with some arbitrary $u_{i^{\ast}}\in \mathbb{R}$, satisfies condition \eqref{equation_n_derivative_positive} for $k=0,1$ and condition \eqref{theorem_conditions_definition_of_mu}. If $0 < x^{\prime}_{i^{\ast}} \leq 1$
   and
   \begin{equation}\label{Proof_lemma_about_stability_of_V_function_condition_for_h}
     0 < h < \frac{1}{3\mu},
   \end{equation}
   where constant $\mu$ is defined in \eqref{theorem_conditions_definition_of_mu}, then for $\bar{A}_{i}, \bar{B}_{i}$ and $x^{\prime}_{i},$ $i = i^{\ast}, i^{\ast}+1, \ldots,$ calculated according to formulas \eqref{SI_def_third}, the inequalities hold true
   \begin{equation}\label{Proof_lemma_about_stavility_of_V_function_A_B_less_then_zero}
     \bar{A}_{i}\frac{u^{2}}{2}+\bar{B}_{i}u \leq 0,\; \forall u\in [0, h],
   \end{equation}
   \begin{equation}\label{Proof_lemma_about_stavility_of_V_function_x'_yends_to_zero}
     0\leq x^{\prime}_{i} \leq \left(\frac{1}{\left(x_{i^{\ast}}^{\prime}\right)^{2}} + \int\limits_{u_{i^{\ast}}}^{u_{i} - h}\mathcal{N}(u)du\right)^{-\frac{1}{2}}.
   \end{equation}
\end{lemma}
\begin{proof}
Let us consider an auxiliary sequence
\begin{equation}\label{Proof_lemma_about_stavility_auxiliary_problem}
  \bar{x}^{\prime}_{i} = \bar{x}^{\prime}_{i-1}\exp\left(s_{1,i}\mathcal{N}(u_{i-1})\left(\bar{x}^{\prime}_{i-1}\right)^{2}h + s_{2,i}\right), \; s_{1, i} \leq -\frac{1}{2},\; s_{2, i} \leq 0,
\end{equation}
$$i= i^{\ast} + 1, i^{\ast} + 2, \ldots; \; \bar{x}^{\prime}_{i^{\ast}} = x^{\prime}_{i^{\ast}}.$$
From the inequality
$$\left(\bar{x}_{i}^{\prime}\right)^{2} - \left(\frac{1}{\left(\bar{x}_{i-1}^{\prime}\right)^{2}} + \mathcal{N}(u_{i-1})h\right)^{-1} = $$
$$ = \left(\bar{x}^{\prime}_{i-1}\right)^{2}\exp\left(2 s_{1,i}\mathcal{N}(u_{i-1})\left(\bar{x}^{\prime}_{i-1}\right)^{2}h + 2 s_{2,i}\right) - \frac{\left(\bar{x}^{\prime}_{i-1}\right)^{2}}{1+\left(\bar{x}^{\prime}_{i-1}\right)^{2}\mathcal{N}(u_{i})h} =$$
$$= \frac{\left(\bar{x}^{\prime}_{i-1}\right)^{2}\left(\exp\left(2 s_{1,i}\mathcal{N}(u_{i-1})\left(\bar{x}^{\prime}_{i-1}\right)^{2}h + 2 s_{2,i}\right)\left(1 +\left(\bar{x}^{\prime}_{i-1}\right)^{2}\mathcal{N}(u_{i-1})h\right) - 1\right)}{1+\left(\bar{x}^{\prime}_{i-1}\right)^{2}\mathcal{N}(u_{i})h} \leq$$
$$\leq \frac{\left(\bar{x}^{\prime}_{i-1}\right)^{2}\left(\exp\left(2 s_{1,i}\mathcal{N}(u_{i-1})\left(\bar{x}^{\prime}_{i-1}\right)^{2}h + 2 s_{2,i}\right)\exp\left(\mathcal{N}(u_{i-1})\left(\bar{x}^{\prime}_{i-1}\right)^{2}h\right) - 1\right)}{1+\left(\bar{x}^{\prime}_{i-1}\right)^{2}\mathcal{N}(u_{i})h} \leq 0$$
it follows that
\begin{equation}\label{Proof_lemma_about_stavility_main_inequality}
  \bar{x}^{\prime}_{i} \leq \left(\frac{1}{\left(\bar{x}_{i-1}^{\prime}\right)^{2}} + \mathcal{N}(u_{i-1})h\right)^{-\frac{1}{2}}.
\end{equation}
Applying inequality \eqref{Proof_lemma_about_stavility_main_inequality} recursively we get the estimate
\begin{equation}\label{Proof_lemma_about_stavility_main_inequality_global}
  \bar{x}^{\prime}_{i} \leq \left(\frac{1}{\left(\bar{x}_{i^{\ast}}^{\prime}\right)^{2}} + \sum\limits_{j = i^{\ast}}^{i-1}\mathcal{N}(u_{j})h\right)^{-\frac{1}{2}} \leq
   \left(\frac{1}{\left(\bar{x}_{i^{\ast}}^{\prime}\right)^{2}} + \int\limits_{u_{i^{\ast}}}^{u_{i-1}}\mathcal{N}(\xi)d\xi\right)^{-\frac{1}{2}} \leq
\end{equation}
$$\leq \left(1 + \int\limits_{u_{i^{\ast}}}^{u_{i}}\mathcal{N}(u)du - \mathcal{N}(u_{i})h\right)^{-\frac{1}{2}}.$$
 To derive the last two inequalities in \eqref{Proof_lemma_about_stavility_main_inequality_global} we exploited the fact that function $\mathcal{N}(u)$ is
 non-decreasing (see \eqref{equation_n_derivative_positive}, $k=1$).
 From \eqref{Proof_lemma_about_stavility_main_inequality_global}, using \eqref{Proof_lemma_about_stability_of_V_function_condition_for_h} we get
 \begin{equation}\label{Proof_lemma_about_stavility_exponent_argument_less_than_1/2}
   0 \leq \mathcal{N}(u_{i})\left(\bar{x}^{\prime}_{i}\right)^{2}h \leq \mathcal{N}(u_{i})h\left(1 + \int\limits_{u_{i^{\ast}}}^{u_{i}}\mathcal{N}(u)du - \mathcal{N}(u_{i})h\right)^{-1} \leq \frac{\mu h}{1-\mu h}\leq \frac{1}{2},
 \end{equation}
 $$i = i^{\ast}, i^{\ast} + 1, \ldots. $$

Inequality \eqref{Proof_lemma_about_stavility_exponent_argument_less_than_1/2} together with \eqref{equation_n_derivative_positive} ($k=0,1$) imply that estimate \eqref{Proof_lemma_about_stavility_main_inequality_global} remains valid
if
\begin{equation}\label{Proof_lemma_about_stavility_substitution_of_s1_s2}
  s_{1,i} = \mathcal{N}(u_{i-1})\left(\bar{x}^{\prime}_{i-1}\right)^{2}h - 1 \leq -\frac{1}{2},\; s_{2,i} = -\frac{h^{2}}{2}\mathcal{N}^{\prime}(u)\left(\bar{x}^{\prime}_{i-1}\right)^{2} \leq 0.
\end{equation}
On the other hand, sequence $\{\bar{x}^{\prime}_{i}\}$ \eqref{Proof_lemma_about_stavility_auxiliary_problem}, with constants $s_{1,i},$ $s_{2,i}$ defined as
\eqref{Proof_lemma_about_stavility_substitution_of_s1_s2}, totally coincide with sequence $\{x^{\prime}_{i}\}$ \eqref{SI_def_third} :
\begin{equation}\label{Proof_lemma_about_stavility_sequence_x_i_expanded_form}
  x_{i}^{\prime} = x_{i-1}^{\prime} \exp\left(-h \mathcal{N}(u_{i-1})\left(x_{i-1}^{\prime}\right)^{2} + \frac{h^{2}}{2}\left(-\mathcal{N}^{\prime}(u_{i-1})\left(x_{i-1}^{\prime}\right)^{2} + 2\left(\mathcal{N}(u_{i-1}) \left(x_{i-1}^{\prime}\right)^{2}\right)^{2}\right)\right) =
\end{equation}
$$= x_{i-1}^{\prime}\exp\left(A_{i-1}\frac{h^{2}}{2} + B_{i-1} h\right) = V(A_{i-1}, B_{i-1}, x^{\prime}_{i-1}, h),\; i=i^{\ast} + 1, i^{\ast} + 2, \ldots.$$
In the light of the latter observation, estimates \eqref{Proof_lemma_about_stavility_of_V_function_x'_yends_to_zero}, immediately follow from \eqref{Proof_lemma_about_stavility_main_inequality_global},
whereas inequalities \eqref{Proof_lemma_about_stavility_of_V_function_A_B_less_then_zero} follow from \eqref{Proof_lemma_about_stavility_substitution_of_s1_s2} and \eqref{SI_def_third}.
\end{proof}

\begin{proof}[Proof of Theorem \ref{Second_theorem_about_approximation_properties_of_SI_method}]
Let the conditions of Theorem \ref{Second_theorem_about_approximation_properties_of_SI_method} are fulfilled. This immediately implies that conditions of Lema \ref{Proof_lemma_about_stability_of_V_function} are fulfilled as well and inequalities \eqref{theorem_two_x_i_prime_relains_below_1} follow from \eqref{Proof_lemma_about_stavility_of_V_function_x'_yends_to_zero} in a trivial way.

Below we implicitly use a fact established in scope of Theorem \ref{Main_theorem_about_approximation_of_IS_method} (whose conditions are fulfilled) that functions $u(x)$ and $u^{\prime}(x)$ are monotonically increasing (see \eqref{monotonicity_of_u} and \eqref{deriv_of_straight_tends_to_infinity}).

At this point we focus on deriving estimates for $|x(u_{i^{\ast}}) - x_{i^{\ast}}|$ and $|x^{\prime}(u_{i^{\ast}}) - x^{\prime}_{i^{\ast}}|.$ This require us to prove some auxiliary inequalities for $x^{\prime}(x)$ and $x^{\prime\prime}(x)$ on $[u(x_{i^{\ast}} - h), u(x_{i^{\ast}} + h)]$ as it follows below.

Using \eqref{theorem_conditions_first_estimation} and \eqref{Theorem_aux_estimete_magnitude} we get the estimate
\begin{equation}\label{Proof_estimation_from_below_for_u'}
  \min\limits_{x\in \left[x_{i^{\ast}} - h, x_{i^{\ast}} + h\right]} u^{\prime}(x) = u^{\prime}\left(x_{i^{\ast}} - h\right) =
  u^{\prime}\left(x_{i^{\ast}} - h\right) - u^{\prime}\left(x_{i^{\ast}}\right) + u^{\prime}\left(x_{i^{\ast}}\right) - u^{\prime}_{i^{\ast}} + u^{\prime}_{i^{\ast}} \geq
\end{equation}
$$ \geq -u^{\prime \prime}(x_{i^{\ast}})h - Ph^{2} + 1 \geq 1 - Ph^{2} -  L_{0}M_{0} h = \tau(h), $$
which immediately implies
\begin{equation}\label{Proof_estimate_for_x_prime}
  \max\limits_{u\in [u(x_{i^{\ast}} - h), u(x_{i^{\ast}} + h)]} x^{\prime}(u) \leq \frac{1}{\tau(h)}.
\end{equation}
At this point, we have to mention that from \eqref{theorem_conditions_first_main_restriction_on_h} and \eqref{Theorem_1_epsilon_definition} it follows that $$\tau(h) > 1-2\varepsilon > 0.$$

Now, imposing a restriction on the magnitude of $h$
\begin{equation}\label{Proof_third_upper_bound_restriction_for_h}
  h \leq \frac{\tau(h)}{P},
\end{equation}
which immediately follows from condition \eqref{theorem_conditions_second_main_restriction_on_h}, we assert that
\begin{equation}\label{Proof_u_i_ast_belongs_to_interval}
  u_{i^{\ast}} \in [u(x_{i^{\ast}} - h), u(x_{i^{\ast}} + h)].
\end{equation}
Indeed, in the light of \eqref{Proof_third_upper_bound_restriction_for_h}, inclusion \eqref{Proof_u_i_ast_belongs_to_interval} can be justified by the two inequalities below
$$u_{i^{\ast}} - u(x_{i^{\ast}} - h) = u_{i^{\ast}} - u(x_{i^{\ast}}) + u(x_{i^{\ast}}) - u(x_{i^{\ast}} - h) \geq h\tau(h) - Ph^{2} \geq 0,$$
$$u(x_{i^{\ast}} + h) - u_{i^{\ast}} = u(x_{i^{\ast}} + h) - u(x_{i^{\ast}}) + u(x_{i^{\ast}}) - u_{i^{\ast}} \geq h\tau(h) - Ph^{2} \geq 0.$$

From \eqref{Proof_u_i_ast_belongs_to_interval} it follows that $u(x_{i^{\ast}}) + \theta (u_{i^{\ast}} - u(x_{i^{\ast}})) \in [u(x_{i^{\ast}} - h), u(x_{i^{\ast}} + h)],$
$\forall \theta \in [0,1].$ With this in mind, and using estimates \eqref{theorem_conditions_first_estimation}, \eqref{Proof_estimate_for_x_prime},  we derive the inequality
\begin{equation}\label{Proof_estimation_x_i_ast}
  |x(u_{i^{\ast}}) - x_{i^{\ast}}| = |x(u_{i^{\ast}}) - x(u(x_{i^{\ast}}))| = x^{\prime}(u(x_{i^{\ast}}) + \theta (u_{i^{\ast}} - u(x_{i^{\ast}})))|u_{i^{\ast}} - u(x_{i^{\ast}})| \leq \frac{Ph^{2}}{\tau(h)}.
\end{equation}
where $0 \leq \theta \leq 1.$

From \eqref{Theorem_aux_estimete_magnitude} and \eqref{Proof_estimate_for_x_prime} it follows that
\begin{equation}\label{Proof_estimation for x''}
  \max\limits_{u\in [u(x_{i^{\ast}} - h), u(x_{i^{\ast}} + h)]} |x^{\prime\prime}(u)| = \max\limits_{u\in [u(x_{i^{\ast}} - h), u(x_{i^{\ast}} + h)]} N(u)u\left(x^{\prime}(u)\right)^{3} \leq
  \frac{L_{0}M_{0}}{\tau(h)^{3}}.
\end{equation}
Inclusion \eqref{Proof_u_i_ast_belongs_to_interval} together with inequalities \eqref{theorem_conditions_first_estimation} and \eqref{Proof_estimation for x''} yield us the estimate
\begin{equation}\label{Proof_estimate for x'}
  |x^{\prime}(u_{i^{\ast}}) - x^{\prime}_{i^{\ast}}| \leq |x^{\prime}(u_{i^{\ast}}) - x^{\prime}(u(x_{i^{\ast}}))| + |x^{\prime}(u(x_{i^{\ast}}))- x^{\prime}_{i^{\ast}}| \leq
\end{equation}
$$\leq \frac{L_{0}M_{0}P}{\tau(h)^{3}}h^{2} + \left|\frac{1}{u^{\prime}(x_{i^{\ast}})} - \frac{1}{u^{\prime}_{i^{\ast}}}\right|
\leq \left(\frac{L_{0}M_{0}}{\tau(h)^{3}} + \frac{1}{1-Ph^{2}}\right)Ph^{2}.$$

At this point, by inequalities \eqref{Proof_estimation_x_i_ast} and \eqref{Proof_estimate for x'} we proved estimates \eqref{theorem_statement_about_x_approximation} and \eqref{theorem_statement_about_x_prime_approximation} respectively for $i = i^{\ast}.$ Below we address the case $i > i^{\ast}.$

Let us now consider a sequence of functions $\{y_{i}(u)\}_{i^{\ast}}^{\infty},$ defined as follows
\begin{eqnarray}
   y_{i}(u) &=& x(u) - V_{i}(u), \; u\in[u_{i}, u_{i+1}],\\ \label{Proof_vector_function_Y_definition}
   V_{i}(u) & = & V\left(\bar{A}_{i}, \bar{B}_{i}, x^{\prime}_{i}, x_{i}, u-u_{i}\right), \nonumber
\end{eqnarray}
where $\bar{A}_{i}, \bar{B}_{i}$ for $i \geq i^{\ast}$ and  $u_{i}, x^{\prime}_{i}, x_{i}$ for $i > i^{\ast}$ are defined according to formulas \eqref{SI_def_third}; to be more specific:
$$   \bar{A}_{i} =  -\mathcal{N}^{\prime}(u_{i})\left(x^{\prime}_{i}\right)^{2} + 2\left(\mathcal{N}(u_{i})\right)^{2}\left(x_{i}^{\prime}\right)^{4},$$
$$   \bar{B}_{i}  =  - \mathcal{N}(u_{i})\left(x_{i}^{\prime}\right)^{2}, i= i^{\ast}, i^{\ast} + 1, \ldots. $$

It is easy to see that $y_{i}(u)$ should satisfy the recurrence system of Cauchy problems
\begin{equation}\label{Proof_equation_fro_Y_i}
  y^{\prime\prime}_{i}(u) = G_{i}(u)y^{\prime}_{i}(u) + F_{i}(u)x^{\prime}(u),
\end{equation}
$$G_{i}(u) = \bar{B}_{i} + \bar{A}_{i}(u - u_{i}), \; F_{i}(u) = - \mathcal{N}(u)\left(x^{\prime}(u)\right)^{2} - \bar{A}_{i}(u-u_{i}) - \bar{B}_{i},\; u\in [u_{i}, u_{i+1}],$$
\begin{equation}\label{Proof_init_condition_for_Y_i}
  y_{i}(u_{i+1}) = y_{i+1}(u_{i+1}), \; y^{\prime}_{i}(u_{i+1}) = y^{\prime}_{i+1}(u_{i+1}), \; u_{i+1} = u_{i} + h, \; i= i^{\ast}, i^{\ast} + 1,\ldots.
\end{equation}
Inequalities \eqref{Proof_estimation_x_i_ast}, \eqref{Proof_estimate for x'} allow us to estimate $|y^{(k)}_{i^{\ast}}(u_{i^{\ast}})|,$ $k=0,1$ in the following way:
\begin{equation}\label{Proof_esqimate_fro_Y_i_ast-1}
  |y_{i^{\ast}}(u_{i^{\ast}})| \leq \frac{Ph^{2}}{\tau(h)},\; |y^{\prime}_{i^{\ast}}(u_{i^{\ast}})| \leq \left(\frac{L_{0}M_{0}}{\tau(h)^{3}} + \frac{1}{1-Ph^{2}}\right)Ph^{2}.
\end{equation}

Using the mean value theorem, we find that
\begin{equation}\label{Proof_F_i_and_mean_value_theorem}
  F_{i}(u) = \frac{\left(u-u_{i}\right)^{2}}{2}\left(-\mathcal{N}^{\prime\prime}(u_{i} + \theta (u-u_{i}))\left(x^{\prime}(u_{i} + \theta (u-u_{i}))\right)^{2}\right. +
\end{equation}
$$+6 \mathcal{N}(u_{i} + \theta (u-u_{i}))\mathcal{N}^{\prime}(u_{i} + \theta (u-u_{i}))\left(x^{\prime}(u_{i} + \theta (u-u_{i}))\right)^{4} -$$
$$\left.-8\left(\mathcal{N}(u_{i} + \theta (u-u_{i}))\right)^{3}\left(x^{\prime}(u_{i} + \theta (u-u_{i}))\right)^{6}\right)
-\mathcal{N}(u_{i})\left(\left(x^{\prime}(u_{i})\right)^{2} - \left(x^{\prime}_{i}\right)^{2}\right) +$$
$$ + (u - u_{i})\left(-\mathcal{N}^{\prime}(u_{i})\left(\left(x^{\prime}(u_{i})\right)^{2} - \left(x^{\prime}_{i}\right)^{2}\right) + 2\left(\mathcal{N}(u_{i})\right)^{2}\left(\left(x^{\prime}(u_{i})\right)^{4} - \left(x^{\prime}_{i}\right)^{4}\right)\right),$$
$$\forall u\in [u_{i}, u_{i+1}],\; i \geq i^{\ast},\; \theta = \theta(u) \in (0,1),$$
which, together with \eqref{Proof_lemma_about_stavility_of_V_function_x'_yends_to_zero} and \eqref{inverse_solution_closed_form}, yields us
an estimate\footnote{Here we use conditions \eqref{equation_n_derivative_positive}, implying that functions $\mathcal{N}(u),$ $\mathcal{N}^{\prime}(u)$ and $\mathcal{N}^{\prime\prime}(u)$ are nondecreasing.}
\begin{equation}\label{Proof_F_i_estimate}
  |F_{i}(u)| \leq\frac{\left(u-u_{i}\right)^{2}}{2}\left(\frac{\mathcal{N}^{\prime\prime}(u)}{\Upsilon_{1}(u-h, h)} + \frac{6\mathcal{N}^{\prime}(u)\mathcal{N}(u)}{\left(\Upsilon_{1}(u-h, h)\right)^{2}} +\frac{8\left(\mathcal{N}(u)\right)^{3}}{\left(\Upsilon_{1}(u-h, h)\right)^{3}}\right) +
\end{equation}
$$ + 2|x^{\prime}(u_{i}) - x^{\prime}_{i}|\left(\frac{\mathcal{N}(u) + h\mathcal{N}^{\prime}(u)}{\left(\Upsilon_{2}(u - 2h)\right)^{\frac{1}{2}}}
+\frac{2h\left(\mathcal{N}(u)\right)^{2}}{\left(\Upsilon_{2}(u -2h)\right)^{\frac{3}{2}}} \right),\; u\in [u_{i}, u_{i+1}],\; i=i^{\ast}, i^{\ast}+1, \ldots,$$
where $\Upsilon_{1}(u, h)$ and $\Upsilon_{2}(u)$ are defined in \eqref{Second_theorem_Upsilon_1_2_definition}.

It is worth mentioning, that to derive \eqref{Proof_F_i_estimate} we estimated $x(u)$ using formula \eqref{inverse_solution_closed_form} and inequality \eqref{Proof_estimate_for_x_prime} as follows:
$$x^{\prime}(u) = \left(\left(x^{\prime}(u_{i^{\ast}})\right)^{-2} + 2\int\limits_{u_{i^{\ast}}}^{u}\mathcal{N}(\xi)d\xi\right)^{-\frac{1}{2}}\leq \left(\tau\left(h\right)^{2} + 2\int\limits_{u_{i^{\ast}}}^{u}\mathcal{N}(\xi)d\xi\right)^{-\frac{1}{2}} = \frac{1}{\sqrt{\Upsilon_{1}(u, h)}}.$$

Solution to IVP \eqref{Proof_equation_fro_Y_i}, \eqref{Proof_init_condition_for_Y_i} can be expressed in the form
\begin{equation}\label{Proof_solution_for_Y_i_prime}
  y_{i}^{\prime}(u) = \int\limits_{u_{i}}^{u}\exp\left(\int\limits_{\xi}^{u}G_{i}(\zeta)d\zeta\right)F_{i}(\xi)x^{\prime}(\xi)d\xi + y_{i-1}^{\prime}(u_{i})\exp\left(\int\limits_{u_{i}}^{u}G_{i}(\zeta)d\zeta\right),
\end{equation}
\begin{equation}\label{Proof_solution_for_Y_i}
  y_{i}(u) = \int\limits_{u_{i}}^{u}y_{i}^{\prime}(\xi)d\xi + y_{i-1}(u_{i}),
\end{equation}
$$u\in [u_{i}, u_{i + 1}], \; i = i^{\ast}, i^{\ast} + 1, \ldots.$$

Using estimates \eqref{Proof_lemma_about_stavility_of_V_function_A_B_less_then_zero} and \eqref{Proof_F_i_estimate}, from \eqref{Proof_solution_for_Y_i_prime} we get the inequalities
$$|y_{i}^{\prime}(u)| \leq \frac{h^{2}}{2}\int\limits_{u_{i}}^{u}\mathcal{T}(\zeta,h)d\zeta +$$ $$+ \frac{h^{2}}{2}\left(\int\limits_{u_{i}}^{u}\Lambda(\zeta, h)d\zeta + 1\right)\sum\limits_{j=i^{\ast}}^{i-1}\prod\limits_{k=j+1}^{i-1}\left(\int\limits_{u_{k}}^{u_{k+1}}\Lambda(\zeta, h)d\zeta + 1\right)\int\limits_{u_{j}}^{u_{j+1}}\mathcal{T}(\zeta,h)d\zeta + $$
$$+ \left(\int\limits_{u_{i}}^{u}\Lambda(\zeta, h)d\zeta + 1\right)\prod\limits_{k=i^{\ast}}^{i-1}\left(\int\limits_{u_{k}}^{u_{k+1}}\Lambda(\zeta, h)d\zeta + 1\right)|y^{\prime}_{i^{\ast}}\left(u_{i^{\ast}}\right)| \leq $$
$$ \leq \left(\int\limits_{u_{i}}^{u}\Lambda(\zeta, h)d\zeta + 1\right)\prod\limits_{k=i^{\ast}}^{i-1}\left(\int\limits_{u_{k}}^{u_{k+1}}\Lambda(\zeta, h)d\zeta + 1\right)\left(|y^{\prime}_{i^{\ast}}\left(u_{i^{\ast}}\right)| + \frac{h^{2}}{2}\int\limits_{u_{i^{\ast}}}^{u}\mathcal{T}(\zeta,h)d\zeta\right)\leq $$
\begin{equation}\label{Proof_estimate_for_y_i_prime}
   \leq \exp\left(\int\limits_{u_{i^{\ast}}}^{u}\Lambda(\zeta, h)d\zeta\right)\left(|y^{\prime}_{i^{\ast}}\left(u_{i^{\ast}}\right)| + \frac{h^{2}}{2}\int\limits_{u_{i^{\ast}}}^{u}\mathcal{T}(\zeta, h)d\zeta\right),\; u\in [u_{i}, u_{i+1}],\; i = i^{\ast}, i^{\ast} + 1, \ldots,
\end{equation}
where $\Lambda(\zeta, h)$ and $\mathcal{T}(\zeta, h)$ are defined in \eqref{Exponenta_term_expression} and \eqref{Second_derivative_term_expression} respectively.
Combining \eqref{Proof_estimate_for_y_i_prime} with \eqref{Proof_solution_for_Y_i}, we get
\begin{equation}\label{Proof_estimate_for_y_i}
|y_{i}(u)| \leq |y_{i^{\ast}}(u_{i^{\ast}})| + \int\limits_{u_{i^{\ast}}}^{u}\exp\left(\int\limits_{u_{i^{\ast}}}^{\eta}\Lambda(\zeta, h)d\zeta\right)\left(|y^{\prime}_{i^{\ast}}\left(u_{i^{\ast}}\right)| + \frac{h^{2}}{2}\int\limits_{u_{i^{\ast}}}^{\eta}\mathcal{T}(\zeta)d\zeta \right)d\eta,
\end{equation}
$$ u\in [u_{i}, u_{i+1}],\; i = i^{\ast}, i^{\ast} + 1, \ldots.$$
Estimates \eqref{theorem_statement_about_x_approximation}, \eqref{theorem_statement_about_x_prime_approximation} follows immediately from \eqref{Proof_esqimate_fro_Y_i_ast-1} and
estimates \eqref{Proof_estimate_for_y_i_prime}, \eqref{Proof_estimate_for_y_i} respectively. The theorem is proved.
\end{proof}

From the SI-method's perspective, Theorems \ref{Main_theorem_about_approximation_of_IS_method} and \ref{Second_theorem_about_approximation_properties_of_SI_method} mean that under certain conditions (mentioned in the theorems) imposed on the IVP \eqref{first_eq}, \eqref{initial_conditions}, the SI method \eqref{SI_def_first}, \eqref{SI_def_second}, \eqref{SI_def_third}, applied to the problem, behaves in a very predictable way: it starts with the "straight" phase \eqref{SI_def_second}, then at some iteration (with index $i^{\ast}$) it switches to the "inverse" phase \eqref{SI_def_third} and remains within the "inverse" phase no matter how many iterations we perform.

\section{Straight-Inverse method for solving BVPs for second order differential equations. }\label{Section_label_SI_method_for_bvp}
\subsection{Preliminary comments}
Introducing the SI-method for solving two-point boundary value problems, we are going to consider the \textit{simple} and \textit{multiple shooting techniques} supplemented with the SI-method for solving IVPs described above.

We avoid discussing the question about the existence of the solution to BVP \eqref{first_eq}, \eqref{boundary_cond}, assuming, for now, that it is granted (i.e. one has ensured that the solution exists before applying the methods proposed below).

 In this section we do not formulate any theoretical statements that guarantee convergence (success) to either of the proposed methods, even if the exact solution to the BVP exists. The question about sufficient conditions for the methods to converge is rather complex and will be addressed in the subsequent publications.

\subsection{SI single shooting method}
Generally speaking, the \textit{simple shooting technique} for solving BVPs can be imagined as a zero-finding algorithm (e.g. simple bisection), applied to some (problem-dependent) function $F(s),$ which can be evaluated by means of an IVP solver (see, for example, \cite[Section 7.3.1]{stoer2002introduction}). Taking into account specifics of the SI IVP solver \eqref{SI_def_first}, \eqref{SI_def_second}, \eqref{SI_def_third}, it might be a bit tricky to construct the corresponding function $F(s)$ for each particular BVP. This is mainly due to the fact that the SI method, even when applied to a one-dimensional problem, operates in a two-dimensional space, treating both $x$ and $u$ as independent variables depending on the situation. Here we would like to illustrate how the function $F(s)$ can be constructed in a simple case when the exact solution $u(x)$ of the BVP \eqref{first_eq}, \eqref{boundary_cond} is monotone and possesses a single {\it boundary layer} near the right end of the interval $[a, b]$ where $u^{\prime}(x) \gg 1.$ Assuming that for $h$ sufficiently small the iterative process \eqref{SI_def_first}, \eqref{SI_def_second}, \eqref{SI_def_third} approximates the solution of IVP \eqref{first_eq}, \eqref{initial_conditions} (in the sense described in Theorems \ref{Main_theorem_about_approximation_of_IS_method}, \ref{Second_theorem_about_approximation_properties_of_SI_method}), we naturally come to a conclusion that if $u^{\prime}_{a}$ is close enough to $u^{\prime}(a),$ then the process will end up in the "inverse" phase \eqref{SI_def_third}. This gives us a key insight on how to define $F(s).$ In the "inverse" phase we do not have a control over the values of $x_{i},$ which makes it practically impossible for us to ensure that $x_{j} = b$ for some $j.$ But we do have a control over the values of $u_{i}$ by means of adjusting $h.$ Modifying the last formula in \eqref{SI_def_third} as follows
\begin{equation}\label{Step_adjustment_formula}
  h_{i}^{\ast} = \sign(x^{\prime}_{i-1}) \left\{
                                            \begin{array}{cc}
                                              \min\{h, |u_{b} - u_{i-1}|\}, & u_{b} - u_{i-1} \neq 0, \\
                                              h, & u_{b} - u_{i-1} = 0, \\
                                            \end{array}
                                          \right.
\end{equation}
 we can guarantee that there exists an index $j=j(x^{\prime}_{0})$ such that $u_{j} = u_{b}$  and the function $F(s)$ can be defined as
 $$F(s) = \left(x_{j}(s) \; | \; u^{\prime}_{0} = s,\; u_{j} = u_{b}\right) - b.$$

\begin{figure}
    \centering
    \begin{tikzpicture}
    \draw[<->] (7,0) node[below]{$x$} -- (0,0) --
    (0,5) node[left]{$u$};
    \draw (0, 0.5) node[left]{$u_{a}$} -- (7, 0.5);
    \draw (0, 0) node[below]{$a$};
    \draw (0, 3.5) node[left]{$u_{b}$} -- (7, 3.5);
    \draw (4.6, -0.1) node[below]{$b$} -- (4.6, 5);
    \draw[very thick] (0,0.5) to [out=0,in=-110] (4.6, 3.5) to [out = 70, in = -98] (5,5);
    \draw (5,5) node[above]{$u(x)$};
    \draw[dashed] (0,0.5) to [out=5,in=-110] (4.0, 3.5) to [out = 80, in = -90] (4.2,5);
    \draw (4.0, 2.7) node[left]{\rotatebox{60}{$x_{0}^{\prime} = s_{1}^{-}$}};
    \draw[dashed] (0,0.5) to [out=7,in=-95] (2.6, 3.5) to [out = 77, in = -95] (2.8,5);
    \draw (2.7, 2.7) node[left]{\rotatebox{70}{$x_{0}^{\prime} = s_{0}^{-}$}};
    \draw[dashed] (0,0.5) to [out=0,in=-110] (5.6, 3.5) to [out = 77, in = -95] (5.9,5);
    \draw (4.75, 2.7) node[right]{\rotatebox{57}{$x_{0}^{\prime} = s_{1}^{+}$}};
    \draw[dashed] (0,0.5) to [out=0,in=-110] (6.7, 3.5) to [out = 77, in = -95] (7.0,5);
    \draw (5.9, 2.7) node[right]{\rotatebox{62}{$x_{0}^{\prime} = s_{0}^{+}$}};
    \draw[-](2.6, 3.5) -- (2.6, 5);
    \draw[<->](2.6, 4.9) -- (4.6, 4.9);
    \draw (3.5, 5) node[above]{$|F(s_{0}^{-})|$};
    \end{tikzpicture}
    \caption{SI single shooting method. Dashed curves describe SI approximations of the solutions to IVP \eqref{first_eq}, \eqref{initial_conditions} for various initial slopes $u_{0}^{\prime}.$ } \label{fig:M2}
    \end{figure}
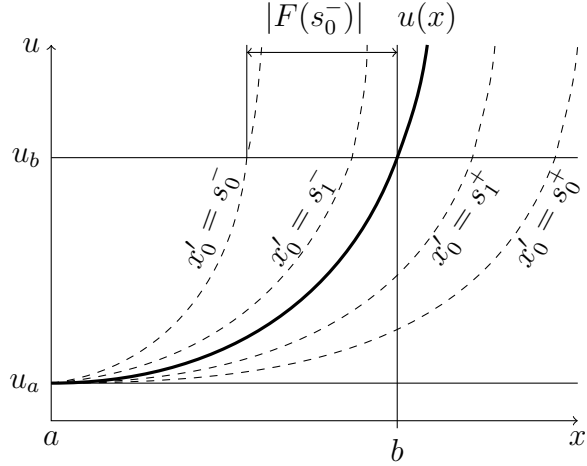

Now if we have $s_{0}^{-}, s_{0}^{+}$ such that $F(s_{0}^{-}) < 0, $ $F(s_{0}^{+}) > 0$
we are ready to run a standard bisection process:
$$s_{i}^{-} = \left\{
              \begin{array}{cc}
                s_{i}, & F(s_{i}) \leq 0, \\
                s_{i-1}^{-} & otherwise, \\
              \end{array}
            \right. \;
            s_{i}^{+} = \left\{
              \begin{array}{cc}
                s_{i}, & F(s_{i}) \geq 0, \\
                s_{i-1}^{+} & otherwise, \\
              \end{array}
            \right. \; s_{i} = \frac{s_{i-1}^{-} + s_{i-1}^{+}}{2},
$$
 which will result in a shrinkage of the distance between $s_{i}^{-}$ and $s_{i}^{+}$ as $i \rightarrow \infty,$  see Fig. \ref{fig:M2} for illustration. We say that the SI {\it single shooting method} converges if
 $$\lim\limits_{i\rightarrow +\infty}F(s_{i}^{-}) = \lim\limits_{i\rightarrow +\infty}F(s_{i}^{+}) = 0.$$
 In the latter case, the truncation of SI mesh \eqref{Definitions_of_omega_mesh}
 $$\Omega(h) = \left\{\left(u^{\prime}_{i}, x^{\prime}_{i}, u_{i}, x_{i}\right), \; i\in\overline{0, j}\; | \; x_{0}=a, u_{0}=u_{a}, x_{j}=b, u_{j}=u_{b}, x_{0}^{\prime} = \lim\limits_{i\rightarrow +\infty}s_{i}^{+}\right\}$$
 is called the SI {\it single shooting approximation} of the solution to BVP \eqref{first_eq}, \eqref{boundary_cond}.

 \begin{proposition}\label{Proposition_about_the_order_of_SI_single_shoothing_approximation}
 If the function $N(u,x)$ and parameter $h$ satisfy conditions of Theorems \ref{Main_theorem_about_approximation_of_IS_method}, \ref{Second_theorem_about_approximation_properties_of_SI_method} and condition \eqref{restrictions_on_initial_conditions} holds true with $u_{a}^{\prime} = \lim\limits_{j\rightarrow +\infty}s_{j}^{+},$ then the SI single shooting approximation of the solution to BVP \eqref{first_eq}, \eqref{boundary_cond} is of order $\mathcal{O}(h^{2})$ in the sense described by the theorems.
 \end{proposition}

\subsection{SI multiple shooting method}
Speaking about the \textit{multiple shooting technique} for solving boundary value problems, we, as a rule, mean a way how the given BVP can be transformed into a system of nonlinear algebraic equations together with an algorithm for solving the system. Below we show how the corresponding system can be constructed using the SI "philosophy".

Assume that we have some initial guess $\Omega_{k}$
\begin{equation}\label{initial_approximation_for_multiple_shooting_method}
\begin{split}
  \Omega_{k} = & \left\{\omega_{k, i} = (u^{\prime}_{k, i}, x_{k,i}^{\prime},  u_{k, i}, x_{k, i}) | \; u_{k,i}^{\prime}\stackrel{def}{=}1/x_{k,i}^{\prime},\; i \in \overline{ 0,N_{k}}; \right.\\
  & \Bigl. x_{k, i} < x_{k, j} \Leftrightarrow i < j < N_{k};\;  u_{k,0} = u_{a}, \; u_{k, N_{k}} = u_{b}, \; x_{k,0} = a, \; x_{k, N_{k}} = b\Bigr\},
\end{split}
\end{equation}
which is a discrete approximation\footnote{The approximation can be constructed using the simple shooting approach described above. } of the exact solution $u(x)$ of the BVP \eqref{first_eq}, \eqref{boundary_cond} in the following sense:
$$ \left.
     \begin{array}{l}
       |u_{k,i}^{\prime}| \leq 1 \Rightarrow\; u(x_{k, i}) \approx u_{k, i},\; u^{\prime}(x_{k, i}) \approx u^{\prime}_{k, i}, \\
       |x_{k,i}^{\prime}| < 1 \Rightarrow x(u_{k,i})\approx x_{k,i}, x^{\prime}(u_{k,i})\approx x^{\prime}_{k,i}, \\
     \end{array}
   \right. \; i \in \overline{ 0,N_{k}}, \; x(\cdot) = u^{-1}(\cdot).
 $$

In what follows we use the notation
\begin{equation}
        h_{k,i} \stackrel{def}{=} x_{k, i + 1} - x_{k, i},\; \bar{h}_{k,i} \stackrel{def}{=} u_{k, i + 1} - u_{k, i},\; i\in \overline{0, N_{k} - 1}
\end{equation}
and require that
\begin{equation}\label{inequality_for_h_bold}
  \max\{h_{k,i}, |\bar{h}_{k,i}|\} \leq h,\; \forall i\in \overline{0, N_{k} - 1}
\end{equation}
for some fixed parameter $h > 0.$

 Combining the general approach, described, for example, in \cite[Section 7.3.5]{stoer2002introduction}, with recurrence formulas \eqref{SI_def_second}, \eqref{SI_def_third}, we transform the initial guess $\Omega_{k}$ into an ordered set of nonlinear equations $$\Gamma_{k} = \{\gamma_{k, i, j}, \; i \in \overline{0, N_{k} - 1},\; j = 0,1\}$$ as it is shown below.

The first two equations can be represented in the form of
\begin{equation}\label{newton_eq_1}
\begin{split}
    & \gamma_{k,0,0}(\omega_{k,0},\; \omega_{k,1}) \stackrel{def}{=}\\
    & \left\{
    \begin{array}{cc}
         U(A(x_{k, 0}, u_{k, 0}, \mathbf{u}^{\prime}_{0}),\;
         B(x_{k, 0}, u_{k, 0}, \mathbf{u}^{\prime}_{0}),\; \mathbf{u}^{\prime}_{0}, u_{k,0}, h_{k,0}) = & \mathbf{u_{1}}, \; |u^{\prime}_{k,0}| \leq 1,\\
         V(\bar{A}(u_{k, 0}, x_{k, 0}, \mathbf{x}^{\prime}_{0}),\;
         \bar{B}(u_{k, 0}, x_{k, 0}, \mathbf{x}^{\prime}_{0}),\; \mathbf{x}^{\prime}_{0}, x_{k,0}, \bar{h}_{k, 0}) = & \mathbf{x_{1}},\; |u^{\prime}_{k,0}| > 1,
    \end{array}\right.\\
\end{split}
\end{equation}
\begin{equation}\label{newton_eq_2}
\begin{split}
    & \gamma_{k,0,1}(\omega_{k,0},\; \omega_{k,1}) \stackrel{def}{=}\\
    & \left\{
    \begin{array}{cc}
         \left.U^{\prime}_{h}(A(x_{k, 0}, u_{k, 0}, \mathbf{u}^{\prime}_{0}),\;
         B(x_{k, 0}, u_{k, 0}, \mathbf{u}^{\prime}_{0}),\; \mathbf{u}^{\prime}_{0}, u_{k,0}, h)\right|_{h = h_{k,0}}, &  |u^{\prime}_{k,0}| \leq 1,\\
         \left.V_{h}^{\prime}(\bar{A}(u_{k, 0}, x_{k, 0}, \mathbf{x}^{\prime}_{0}),\;
         \bar{B}(u_{k, 0}, x_{k, 0}, \mathbf{x}^{\prime}_{0}),\; \mathbf{x}^{\prime}_{0}, x_{k,0}, h)\right|_{h = \bar{h}_{k, 0}}, & \; |u^{\prime}_{k,0}| > 1,
    \end{array}\right. = \\
    & = \left\{\begin{array}{cc}
        \mathbf{u}^{\prime}_{1}, & |u^{\prime}_{k, 0}| \leq 1,\; |u^{\prime}_{k, 1}| \leq 1, \\
        1/\mathbf{x}^{\prime}_{1}, & |u^{\prime}_{k, 0}| \leq 1,\; |u^{\prime}_{k, 1}| > 1,\\
        \mathbf{x}^{\prime}_{1}, & |u^{\prime}_{k, 0}| > 1,\; |u^{\prime}_{k, 1}| > 1,\\
        1/\mathbf{u}^{\prime}_{1}, & |u^{\prime}_{k, 0}| > 1,\; |u^{\prime}_{k, 1}| \leq 1.\\
    \end{array}\right.
\end{split}
\end{equation}
The rest $2(N_{k} - 1)$ equations follow the general pattern described below:
\begin{equation}\label{newton_eq_3}
\begin{split}
    & \gamma_{k,i,0}(\omega_{k,i - 1},\; \omega_{k,i},\; \omega_{k,i + 1}) \stackrel{def}{=}\\
    & \left\{
    \begin{array}{ll}
         U(A(x_{k, i}, \mathbf{u}_{i}, \mathbf{u}^{\prime}_{i}),\;
         B(x_{k, i}, \mathbf{u}_{i}, \mathbf{u}^{\prime}_{i}),\; \mathbf{u}^{\prime}_{i}, \mathbf{u}_{i}, h_{k,i}), & |u^{\prime}_{k,i}| \leq 1,\; |u^{\prime}_{k,i - 1}| \leq 1,\\
         U(A(\mathbf{x}_{i}, u_{k,i}, \mathbf{u}^{\prime}_{i}),\;
         B(\mathbf{x}_{i}, u_{k,i}, \mathbf{u}^{\prime}_{i}),\; \mathbf{u}^{\prime}_{i}, u_{k, i}, x_{k,i + 1} - \mathbf{x}_{i}), & |u^{\prime}_{k,i}| \leq 1,\; |u^{\prime}_{k,i - 1}| > 1,\\
         V(\bar{A}(u_{k, i}, \mathbf{x}_{i}, \mathbf{x}^{\prime}_{i}),\;
         \bar{B}(u_{k, i}, \mathbf{x}_{i}, \mathbf{x}^{\prime}_{i}),\; \mathbf{x}^{\prime}_{i}, \mathbf{x}_{i}, \bar{h}_{k, i}), & |u^{\prime}_{k,i}| > 1,\; |u^{\prime}_{k,i - 1}| > 1,\\
         V(\bar{A}(\mathbf{u}_{i}, x_{k,i}, \mathbf{x}^{\prime}_{i}),\;
         \bar{B}(\mathbf{u}_{i}, x_{k, i}, \mathbf{x}^{\prime}_{i}),\; \mathbf{x}^{\prime}_{i}, x_{k, i}, u_{k, i + 1}- \mathbf{u}_{i}), & |u^{\prime}_{k,i}| > 1,\; |u^{\prime}_{k,i - 1}| \leq 1,\\
    \end{array}\right. =\\
    &= \left\{\begin{array}{cc}
        \left\{\begin{array}{cc}
        \mathbf{u}_{i+1}, & i \in \overline{1, N_{k} - 2} \\
        u_{k, N_{k}} = b, & i = N_{k}-1,
        \end{array}\right., & |u^{\prime}_{k,i}| \leq 1, \\
        \left\{\begin{array}{cc}
            \mathbf{x}_{i+1}, & i \in \overline{1, N_{k} - 2} \\
             x_{k, N_{k}} = b, & i = N_{k}-1,
        \end{array}\right. & |u^{\prime}_{k,i}| > 1,\\
    \end{array}
    \right.
\end{split}
\end{equation}
\begin{equation}\label{newton_eq_4}
\begin{split}
    & \gamma_{k,i,1}(\omega_{k,i - 1},\; \omega_{k,i},\; \omega_{k,i + 1}) \stackrel{def}{=}\\
    & \left\{
    \begin{array}{ll}
         \left.U^{\prime}_{h}(A(x_{k, i}, \mathbf{u}_{i}, \mathbf{u}^{\prime}_{i}),\;
         B(x_{k, i}, \mathbf{u}_{i}, \mathbf{u}^{\prime}_{i}),\; \mathbf{u}^{\prime}_{i}, \mathbf{u}_{i}, h)\right|_{h = h_{k,i}}, & |u^{\prime}_{k,i - 1}| \leq 1,\;  |u^{\prime}_{k,i}| \leq 1,\\
         \left.U^{\prime}_{h}(A(\mathbf{x}_{i}, u_{k, i}, \mathbf{u}^{\prime}_{i}),\;
         B(\mathbf{x}_{i}, u_{k, i}, \mathbf{u}^{\prime}_{i}),\; \mathbf{u}^{\prime}_{i}, u_{k,i}, h)\right|_{h = x_{k,i + 1} - \mathbf{x}_{i}}, &  |u^{\prime}_{k,i - 1}| > 1,|u^{\prime}_{k,i}| \leq 1,\\
         \left.V_{h}^{\prime}(\bar{A}(u_{k, i}, \mathbf{x}_{i}, \mathbf{x}^{\prime}_{i}),\;
         \bar{B}(u_{k, i}, \mathbf{x}_{i}, \mathbf{x}^{\prime}_{i}),\; \mathbf{x}^{\prime}_{i}, \mathbf{x}_{i}, h)\right|_{h = \bar{h}_{k, i}}, & |u^{\prime}_{k,i - 1}| > 1, \; |u^{\prime}_{k,i}| > 1,\\
         \left.V_{h}^{\prime}(\bar{A}(\mathbf{u}_{i}, x_{k, i}, \mathbf{x}^{\prime}_{i}),\;
         \bar{B}(\mathbf{u}_{i}, x_{k, i}, \mathbf{x}^{\prime}_{i}),\; \mathbf{x}^{\prime}_{i}, x_{k,i}, h)\right|_{h = u_{k, i + 1} - \mathbf{u}_{i}}, & |u^{\prime}_{k,i - 1}| \leq 1, \; |u^{\prime}_{k ,i}| > 1,
    \end{array}\right. = \\
    & = \left\{\begin{array}{cc}
        \mathbf{u}^{\prime}_{i + 1}, & |u^{\prime}_{k, i}| \leq 1,\; |u^{\prime}_{k, i + 1}| \leq 1, \\
        1/\mathbf{x}^{\prime}_{i + 1}, & |u^{\prime}_{k, i}| \leq 1,\; |u^{\prime}_{k, i + 1}| > 1,\\
        \mathbf{x}^{\prime}_{i + 1}, & |u^{\prime}_{k, i}| > 1,\; |u^{\prime}_{k, i + 1}| > 1,\\
        1/\mathbf{u}^{\prime}_{i + 1}, & |u^{\prime}_{k, i}| > 1,\; |u^{\prime}_{k, i + 1}| \leq 1,\\
    \end{array}\right.
\end{split}
\end{equation}
$i \in \overline{1, N_{k} - 1},$ where
\begin{equation}
\begin{split}
    A(x, u, u^{\prime}) = & N^{\prime}_{u}(u, x)u^{\prime} + N_{x}^{\prime}(u, x),\\
    B(x, u, u^{\prime}) = & N(u, x),\\
    \bar{A}(u, x, x^{\prime}) = & - \left((N^{\prime}_{u}(u, x) + N_{x}^{\prime}(u, x)x^{\prime})u +  N(u, x)\right)\left(x^{\prime}\right)^{2} + 2 \left(N(u, x)u\right)^{2}\left(x^{\prime}\right)^{4}, \\
    \bar{B}(u, x, x^{\prime}) = & - N(u, x)u\left(x^{\prime}\right)^{2},\\
\end{split}
\end{equation}
the variables in bold describe unknowns and $$\mathbf{u}^{\prime}_{i} \stackrel{def}{=} 1/\mathbf{x}^{\prime}_{i}, \; \forall i \in \overline{0, N_{k}}.$$

As one can see, the equations are dependent on the absolute values of $u^{\prime}_{k,i},$ which, according to our assumption about the approximation properties of $\Omega_{k},$ characterize rapidity of variation of the unknown solution $u(x)$ at different points of segment $[a, b].$ This follows the general idea of the \textit{straight-inverse} approach, consisting in switching between the \textit{straight} (i.e. $u(x)$) and \textit{inverse} (i.e. $x(u)$) solutions depending on which of the two behaves better (that is, possesses lower variation in a vicinity of a given point).

Applying a single iteration of the generalized Newton's method (see, for example, \cite[p. 293]{stoer2002introduction}) to the system $\Gamma_{k}$ \eqref{newton_eq_1}, \eqref{newton_eq_2}, \eqref{newton_eq_3}, \eqref{newton_eq_4} we get a new set $\Omega_{k + 1}$ as a combination of $\Omega_{k}$ and the results brought by the Newton's method iteration, assuming that
$$u_{k+1, i} \stackrel{def}{=} \mathbf{u}_{i}^{(1)},\; u^{\prime}_{k+1, i} \stackrel{def}{=} \mathbf{u}^{\prime(1)}_{i},\;x_{k+1, i} \stackrel{def}{=} \mathbf{x}_{i}^{(1)},\; x^{\prime}_{k+1, i} \stackrel{def}{=} \mathbf{x}^{\prime(1)}_{i}$$
wherever it is relevant. Here the variables in bold with superscript denote the first approximation of the Newton's method applied to the system $\Gamma_{k}.$ In practice, it may happen that the set $\Omega_{k + 1},$ obtained in such a way, needs to be sorted out (to fulfill the requirement $ x_{k+1, i} < x_{k+1, j} \Leftrightarrow i < j < N_{k+1}$) and then refined by the linear interpolation (in order to ensure inequality \eqref{inequality_for_h_bold} for $k$ incremented). Once this is done, approximation $\Omega_{k + 1}$ can be used to construct a new system $\Gamma_{k+1}$ which, after applying another iteration of the Newton's method to it, yields us $\Omega_{k + 2}$ and so on and so forth. If the process can be continued for an arbitrary number of iterations (i.e. the corresponding Jacobian matrices, needed to execute the Newton's iterations, are all nonsingular) and
$$\lim\limits_{j \rightarrow +\infty}\|\Omega_{j} - \Omega_{j+1}\| \stackrel{def}{=} \lim\limits_{j \rightarrow +\infty}\sum_{i} \|\omega_{j+1, i} - \omega_{j,i}\| = 0,$$
then we say that the SI {\it multiple shooting method} is convergent and the limiting mesh $\lim\limits_{j \rightarrow +\infty}\Omega_{j}$ is said to be the SI {\it multiple shooting approximation} of the solution to BVP \eqref{first_eq}, \eqref{boundary_cond}.
\begin{proposition}\label{Proposition_about_the_order_of_approximation_of_SI_multiple_shooting_method}
  If the function $N(u, x)$ and parameter $h$ satisfy conditions of Theorems \ref{Main_theorem_about_approximation_of_IS_method}, \ref{Second_theorem_about_approximation_properties_of_SI_method}
  and condition \eqref{restrictions_on_initial_conditions} holds true with $u_{a}^{\prime} = \lim\limits_{j \rightarrow +\infty}x_{j,0}^{\prime}$ then the SI multiple shooting approximation
  of the solution to BVP \eqref{first_eq}, \eqref{boundary_cond} is of order $\mathcal{O}(h^{2})$ in the sense described by the theorems.
\end{proposition}

\section{Numerical examples}\label{Section_label_numerical_examples}

\subsection{Initial value problem}

In the current sub-section we examine the SI-method \eqref{SI_def_first}, \eqref{SI_def_second}, \eqref{SI_def_third} for solving IVPs by applying it to the Cauchy problem \eqref{troesch_eq}, \eqref{initial_conditions} with $a = 0,$ $u_{a} = 0,$ $u_{a}^{\prime} = 0.1.$ One can easily ensure that the problem satisfies conditions of Theorems \ref{Main_theorem_about_approximation_of_IS_method} and \ref{Second_theorem_about_approximation_properties_of_SI_method}. Numerical results corresponding to different values of  $\lambda$ and $h$ are presented in Tab. \ref{Table_IVP_results_theorem_1}, \ref{Table_IVP_results_theorem_2}.

\begin{table}[!htbp]
\footnotesize
\centering
\begin{tabular}{||c c c c c c c c||}
 \hline
$\lambda$ & $h$ & $i^{\ast}$ & $x_{i^{\ast}}$ & $u_{i^{\ast}}$ & $u_{i^{\ast}}^{\prime}$ & $|u_{i^{\ast}} - u(x_{i^{\ast}})|$ & $|u^{\prime}_{i^{\ast}} - u^{\prime}(x_{i^{\ast}})|$ \\ [0.5ex]
 \hline\hline
  & 1e-1 & 15 & 1.5 & 0.5108552223 & 1.0700488967 & 2.6e-4 & 1.9e-3 \\
 2 & 1e-2 & 147 & 1.47 & 0.4800085101 & 1.0022994311 & 2.5e-6 & 1.7e-5 \\
  & 1e-3 & 1469 & 1.469 & 0.4790098303 & 1.0000906016 & 2.6e-8 & 1.7e-7 \\
  & 1e-4 & 14690 & 1.469 & 0.4790098559 & 1.0000907722 & 1.4e-10 & 1.3e-9 \\
 \hline\hline
  & 1e-2 & 37 & 0.37 & 0.1225264682 & 1.0246219988 & 1.0e-5 & 2.9e-4 \\
 8 & 1e-3 & 368 & 0.368 & 0.1205049349 & 1.0067836140 & 1.0e-7 & 2.8e-6 \\
  & 1e-4 & 3673 & 0.3673 & 0.1198024787 & 1.0005354415 & 9.9e-10 & 2.8e-8 \\
 \hline
\end{tabular}
\caption{Approximation errors of the SI-method, applied to IVP \eqref{troesch_eq}, \eqref{initial_conditions}, which correspond to different values of $h.$ The errors are calculated at the last point, $x_{i^{\ast}},$ of the "straight" phase of the SI method (see Theorem \ref{Main_theorem_about_approximation_of_IS_method}). The reference values $u(x_{i^{\ast}})$ and $u^{\prime}(x_{i^{\ast}})$ are calculated using the {\it dverk78} algorithm implemented in Maple 2016.  }
\label{Table_IVP_results_theorem_1}
\end{table}

\begin{table}[!htbp]
\footnotesize
\centering
\begin{tabular}{||c c c c c c c c||}
 \hline
$\lambda$ & $h$ & $i$ & $x_{i}$ & $u_{i}$ & $x^{\prime}_{i}$ & $|x_{i} - x(u_{i})|$ & $|x^{\prime}_{i} - x^{\prime}(u_{i})|$ \\ [0.5ex]
 \hline\hline
  & 1e-1 & 20 & 1.8072353083 & 1.0 & 0.4262211108 & 1.0e-3 & 1.1e-3 \\
 2 & 1e-2 & 199 & 1.8062219401 & 1.0 & 0.4250841708 & 1.1e-5 & 9.7e-6 \\
  & 1e-3 & 1990 & 1.8062111449 & 1.0 & 0.4250746074 & 1.1e-7 & 9.5e-8 \\
  & 1e-4 & 19900 & 1.8062110370 & 1.0 & 0.4250745138 & 1.1e-9 & 9.3e-10 \\
 \hline\hline
  & 1e-2 & 125 & 0.5434971101 & 1.0 & 1.832181142e-2 & 5.9e-5 & 5.7e-8 \\
 8 & 1e-3 & 1248 & 0.5434390645 & 1.0 & 1.832175495e-2 & 5.8e-7 & 5.4e-10 \\
  & 1e-4 & 12475 & 0.5434384906 & 1.0 & 1.8321754416e-2 & 5.9e-9 & 5.0e-12 \\
 \hline
\end{tabular}
\caption{Approximation errors of the SI-method, applied to IVP \eqref{troesch_eq}, \eqref{initial_conditions}, which correspond to different values of $h.$ The errors are calculated at the point $u_{i} = 1.0$ of the "inverse" phase of the method (see Theorem \ref{Second_theorem_about_approximation_properties_of_SI_method}). The reference values $x(u_{i})$ and $x^{\prime}(u_{i})$ are calculated via formula \eqref{inverse_solution_closed_form} within Maple 2016 environment.}
\label{Table_IVP_results_theorem_2}
\end{table}

The data presented in Tab. \ref{Table_IVP_results_theorem_1} confirms the predictions of Theorem \ref{Main_theorem_about_approximation_of_IS_method} about the order of approximation of the SI method during its "straight" phase (which corresponds to interval $[0, x_{i^{\ast}}]$). To get a better understanding of how precise the error estimates of Theorem \ref{Main_theorem_about_approximation_of_IS_method} are, we would like to evaluate them for the case of $\lambda = 2.$
As it was pointed out in \cite{ROBERTS1976291_Troesch_closed_form_solution}, the initial value problem associated with \eqref{troesch_eq} has a pole approximately in
$$  x_{\infty} = \frac{1}{\lambda}\ln\left(\frac{8}{u^{\prime}(0)}\right).$$
This allows us to get an approximation for $S^{\ast}$
$$  S^{\ast} \approx x_{\infty} = 0.5\ln\left(80\right) \approx 2.191013318.$$
At the same time, Remark \ref{Remark_about_estimate_for_S} allows us to lower the value of $S^{\ast},$ taking into account that the right end of the interval of interest, $x_{i^{\ast}},$ does not exceed $1.5:$
  $$S^{\ast} = \min\{1.5, 2.191013318\} = 1.5.$$
Now using Remark \ref{Remark_to_main_theorem} and taking into account that in case of problem \eqref{troesch_eq}, \eqref{initial_conditions}
$$\Phi(u) = \cosh(\lambda u) - \cosh(\lambda u_{a}),$$ we can calculate $M_{0}$ via the formula
  $$M_{0} = \frac{1}{\lambda}\cosh^{-1}\left(\frac{1}{2}\left((1+3\varepsilon)^{2} - \left(u^{\prime}_{a}\right)^{2}\right) + \cosh(\lambda u_{a})\right).$$
Assuming that $$\varepsilon = 0.01,$$ we get
  $$M_{0} \approx \frac{1}{2}\cosh^{-1}\left(\frac{1}{2}\left(1.3^{2} - 0.1^{2}\right) + \cosh(0.2)\right) \approx 0.5010350625.$$
With the value of $M_{0}$ available, we are in the position to evaluate $ L_{i}$ via formulas \eqref{Proof_L_constant_definition}:
$$   L_{0} = \max\limits_{|u| < M_{0} + \varepsilon}\left|\frac{\lambda \sinh(\lambda u)}{u}\right| \approx 4.733711073,$$

$$   L_{1} = \max\limits_{|u| < M_{0} + \varepsilon}\left|\left(\frac{\lambda \sinh(\lambda u)}{u}\right)^{\prime}\right| \approx 3.021065783,$$

$$   L_{2} = \max\limits_{|u| < M_{0} + \varepsilon}\left|\left(\frac{\lambda \sinh(\lambda u)}{u}\right)^{\prime\prime}\right| \approx 7.11152335,$$

Finally, using formula \eqref{Remark_1_re-definition_of_P_constant}, we get

\begin{equation}\label{Example_P_values_calculated}
  P(h=10^{-4}) \approx 7036.8.
\end{equation}

According to Theorem \ref{Main_theorem_about_approximation_of_IS_method}, value $P$ \eqref{Example_P_values_calculated} gives us an error estimate of the SI method
on the interval $[0, x_{i^{\ast}}]$ for $h=10^{-4},$ see \eqref{theorem_conditions_first_estimation}. Turning back to the data from Tab. \ref{Table_IVP_results_theorem_1}, one can conclude that the error estimates of the theorem are much higher than they potentially can be.

A conclusion similar to the one above can be made when evaluating error estimates of Theorem \ref{Second_theorem_about_approximation_properties_of_SI_method} and comparing them to the corresponding error values from Tab. \ref{Table_IVP_results_theorem_2}. At the same time, the predictions of the theorem about the order of approximation with respect to $h$ are in perfect coherence with the numerical data.

It is worth mentioning, that to calculate numerical data presented in Tab. \ref{Table_IVP_results_theorem_2} we used recurrence equalities \eqref{SI_def_third} enhanced by the step adjustment formula \eqref{Step_adjustment_formula} with $u_{b} = 1.$ This allowed us to achieve an absolute equality $u_{i} = 1$ for the corresponding $i$ from Tab. \ref{Table_IVP_results_theorem_2}.

\subsection{Boundary value problem}

Below we present and discuss numerical results of the SI-method applied to the Troesch's problem \eqref{troesch_eq}, \eqref{troesch_bound_cond}.

Both, the SI single and multiple shooting methods demonstrate convergence when applied to the Troesch's problem. As one might expect, the state of convergence and its rate become more and more dependent on the "quality" of the initial guess as $\lambda$ increases. In practice, good results, in terms of efficiency, are obtained when using a combination of the two methods so that a few iterations of the single shooting method (whose region of convergence is not that sensitive to the magnitude of $\lambda$ but the rate of convergence is quite moderate) provide an initial guess \eqref{initial_approximation_for_multiple_shooting_method} for the multiple shooting algorithm (which possesses a rather high convergence rate provided that the initial guess is precise enough). This combination was used to calculate the numerical results presented below. At the same time, for the sake of analysis, it is quite safe to assume that the results are calculated by the SI single shooting method alone: in terms of accuracy the difference is negligible.

Initial slopes $u^{\prime}(0)$ corresponding to different values of $\lambda$ and calculated by different methods are presented in Tab. \ref{Table_initial_slopes}. The two rightmost columns of the table contain the slopes calculated by the SI-method with different values of step size $h.$ Comparing the results of the SI-method to those calculated by the other methods, we see that the order of approximation of the SI-method with respect to $h$ is very close to 2, which is coherent with Propositions \ref{Proposition_about_the_order_of_SI_single_shoothing_approximation} and \ref{Proposition_about_the_order_of_approximation_of_SI_multiple_shooting_method}.

\begin{table}[!htbp]
\footnotesize
\centering
\begin{tabular}{||c c c c c c||}
 \hline
$\lambda$ & \cite{Gen_sol_of_TP} & \cite{Kutniv_Makarov_IVP_solvers} & Maple 2016 \footnote{Using numeric "dsolve" procedure with "abserr = 1e-12"} & SI-method, $h = 10^{-4}$ & SI-method, $h = 10^{-5}$ \\ [0.5ex]
 \hline\hline
 2  &  0.5186322404 & -- & 0.518621219269 & 0.518621219577035 & 0.518621219272419  \\
 3  &  0.255607567 & -- & 0.255604215562 & 0.255604216455332 & 0.255604215571849 \\
 5  &  4.575046433e-02 & -- & 4.575046140632e-02 & 4.575046196263e-02 & 4.575046141188e-02 \\
 8  &  2.587169418e-03 & -- & 2.587169418963e-3 & 2.587169500425e-03 & 2.587169419777e-03 \\
 20 &  1.648773182e-08 & 1.6487734e-8 & -- & 1.648773647e-08 & 1.648773188e-00 \\
 30 &  7.486093793e-13 & 7.4861194e-13 & -- & 7.486098431e-13 & 7.486093844e-13 \\
 50 &  1.542999878e-21 & 1.5430022e-21 & -- & 1.543002448e-21 & 1.542999906e-21 \\
 61 &  -- & 2.5770722e-26 & -- & 2.577078525e-26 & 2.577072299e-26 \\
 100 & 2.976060781e-43 & -- & -- & 2.976075557e-043 & 2.976060927e-043 \\
 \hline
\end{tabular}
\caption{Values of $u^{\prime}(0)$ for the Troesch's problem calculated by different approaches.}
\label{Table_initial_slopes}
\end{table}

\begin{table}[!htbp]
\footnotesize
\centering
\begin{tabular}{||c c c c c c||}
 \hline
$\lambda$ & \cite{Gen_sol_of_TP} & Other & Maple 2016 \footnote{Using numeric "dsolve" procedure with "abserr = 1e-12"} & SI-method, $h = 10^{-4}$ & SI-method, $h = 10^{-5}$ \\ [0.5ex]
 \hline\hline
 2 & 2.406790318 & 2.406939711 \cite{ROBERTS1976291_Troesch_closed_form_solution} & 2.406939831247 & 2.40693982969129 & 2.4069398312315 \\
 3 & 4.266151411 & 4.266222862 \cite{Chang20103303} & 4.266222861803 & 4.26622285457896 & 4.2662228617306 \\
 5 & 12.10049478 & 1.210049546 \cite{Chang20103303} & 12.1004954507778 & 12.1004954359128 &  12.1004954506293 \\
 8 & 54.57983465 & 5.457983447 \cite{Chang20103303} & 54.5798344555735 &  54.5798344412402 &  54.5798344554302  \\
 10 & 148.4064126 & 148.4064212 \cite{Chang20103303} & -- & 148.406421145524  & 148.406421155906 \\
 20 & 22026.29966 & 22026.4657 \cite{Snyman_Troesch_problem_solution} & -- & 22026.4657494062 & 22026.4657494068 \\
 30 & -- & -- & -- &  3269017.37247181 & 3269017.3724718 \\
 50 & -- & -- & -- & 72004899337.3858 & 72004899337.386 \\
 \hline
\end{tabular}
\caption{Values of $u^{\prime}(1)$ for the Troesch's problem calculated by different approaches.}
\label{Table_final_slopes}
\end{table}

The order of the SI-method's error with respect to $h$ near the right boundary point can be estimated empirically from Tab. \ref{Table_final_slopes} which contains values of
$u^{\prime}(1)$ calculated by different methods for different values of $\lambda.$ The two rightmost columns of the table contain values of $u^{\prime}(1)$
calculated by the SI-method with different values of step size $h.$ Examining the table, we should keep in mind that the values calculated by other (than SI) methods are actually inverse to those approximated by the SI-method, i.e. on the segment where derivative of the unknown function $u(x)$ gets bigger than $1$ the method approximates
values of $x^{\prime}(u) = \frac{1}{u^{\prime}(x)}.$ Nevertheless, we still can see that the method's error is of order $2$ with respect to $h,$ just as it is predicted by Propositions \ref{Proposition_about_the_order_of_SI_single_shoothing_approximation} and \ref{Proposition_about_the_order_of_approximation_of_SI_multiple_shooting_method}.

\begin{table}[!htbp]
\footnotesize
\centering
\begin{tabular}{||c c c c c||}
 \hline
 Value & \cite{Chang20103043} & \cite{Gen_sol_of_TP} & SI-method, $h = 10^{-4}$ & SI-method, $h = 10^{-5}$ \\ [0.5ex]
 \hline\hline
 $u(0.1)$ & 4.211183679705e-05 & 4.211189927237e-05 & 4.21119023173e-05 & 4.21118993037e-05  \\
 $u(0.2)$ & 1.299639238293e-04 & 1.299641158237e-04 & 1.29964125220e-04 & 1.29964115920e-04  \\
 $u(0.3)$ & 3.589778855481e-04 & 3.589784013896e-04 & 3.58978427345e-04 & 3.58978401657e-04  \\
 $u(0.4)$ & 9.779014227050e-04 & 9.779027718029e-04 & 9.77902842508e-04 & 9.77902772532e-04  \\
 $u(0.5)$ & 2.659017178062e-03 & 2.659020490351e-03 & 2.659020682593-03 & 2.65902049234e-03 \\
 $u(0.999)$ & 8.889931171768e-01 & 8.889931181558e-01 & 8.89035025083e-01\footnote{For $x = 0.999000491899$} & 8.88994612232e-01\footnote{For $x = 0.999000017539$}  \\
 \hline
\end{tabular}
\caption{Solution to the Troesch's problem with $\lambda = 10$ evaluated at multiple points inside interval $(0, 1)$ via different approaches.}
\label{table_comparison_to_intermediate_values}
\end{table}

Tab. \ref{table_comparison_to_intermediate_values} presents approximations of the solution $u(x)$ to the Troesch's problem
calculated by different methods at points other than the end points of interval $[0, 1].$ Comparing the results obtained by the SI-method for different
values of $h$ with those obtained by other methods, we can conclude that the order of the SI-method's error with respect to $h$ is still very close $2.$ It is worth mentioning, that because of specifics of the SI-method,
one cannot have a control over the points $x_{i}$ belonging to the rightmost part of the interval $[0, 1],$ where the absolute value of the derivative $u^{\prime}(x)$
exceeds $1.$ In the latter case, the method "works" with the inverse function $x(u)$ and it is rather possible to choose points $u_{i}$
where to calculate the approximation of $x(u) = u^{-1}(u)$. This explains why the bottom row in Tab.
\ref{table_comparison_to_intermediate_values} contains approximations by the SI-method for value $x$ close but not equal to $0.999.$

\begin{table}[!htbp]
\footnotesize
\centering
\begin{tabular}{||c c c c c ||}
 \hline
 Source & $u^{\prime}(0)$ & $\|\Omega(h)\|\footnote{Number of knots in the final mesh.}$ & CPU time, sec. & Rel. diff. to \cite{Gen_sol_of_TP}\footnote{Relative difference as compared to $u^{\prime}(0)$ calculated in \cite{Gen_sol_of_TP}. }  \\ [0.5ex]
 \hline\hline
 SI-method, $h = 10^{-2}$ & 3.141990565e-43 & 240 & 0.022 & 5.6e-2 \\
 SI-method, $h = 10^{-3}$ & 2.977378936e-43 & 2208 & 0.054 & 4.4e-4 \\
 SI-method, $h = 10^{-4}$ & 2.976075557e-43 & 21753 & 0.275 & 5.0e-6 \\
 SI-method, $h = 10^{-5}$ & 2.976060927e-43 & 203143 & 2.135 & 4.9e-8 \\
 SI-method, $h = 10^{-6}$ & 2.976060782e-43 & 2081478 & 16.05 & 3.4e-10 \\
 \cite{Gen_sol_of_TP} & 2.976060781e-43 & -- & -- & 0.0 \\
 \hline
\end{tabular}
\caption{Solution to the Troesch's problem with $\lambda = 100$.}
\label{table_performance_results}
\end{table}

Tab. \ref{table_performance_results} allows us to get an insight about the performance of the SI-method and its complexity. The absolute values of execution time listed in the table are obtained on a laptop with CPU Intel(R) Core(TM) i3-3120M, 2.5 GHz
and 8 Gb of RAM, using the single thread implementation available at \url{https://github.com/imathsoft/MathSoftDevelopment} . The dependency between the execution time and the number of knots seems to be close to a linear one, which gives us an evidence that the complexity of the algorithm can be characterized as $\mathcal{O}(\|\Omega(h)\|).$
On the other hand, it is easy to notice that the dependency between
the number of knots, $\|\Omega(h)\|,$ and the step size $h,$ for the Troesch's problem, can be described by the approximate equality $\|\Omega(h)\| \approx 2/h.$ The latter observation allows us
 to estimate complexity of the algorithm applied to the Troesch's problem as $\mathcal{O}(2/h)$.
More thorough investigation of the SI-method's complexity remains beyond the scope of the present paper and is left to the subsequent publications. Potentially, the implementation of the SI-method can be speeded up by parallelization of some subroutines.

\section{Conclusions}\label{Section_label_conclusions}
The SI-method presented in the paper can be considered as a particular implementation of a quite general idea about switching between "straight" and "inverse" problems when
one of them becomes essentially more difficult in terms of numerical calculations than the other one. The approach presented here can be quite
easily modified and applied to ordinary differential equations of different types, by choosing different step functions $U(s)$ and $V(s).$

The particular version of the SI-method presented above, is quite straightforward and efficient in terms of programming. One of its possible c++ implementations is available
at GitHub \footnote{https://github.com/imathsoft/MathSoftDevelopment } and can be used for solving problems other than the Troesch's problem exploited in the present paper.

The results of numerical examples, based on the Troesch's problem, clearly show that the proposed implementation of the SI-method behaves very well, as compared
to the other approaches, in terms of both accuracy and efficiency. It is worth mentioning that this is despite the fact
that the SI-method is general and does not have anything in it which is designed specifically for the purpose of solving the Troesch's problem
(as it is in some other approaches referenced in Section \ref{Section_label_numerical_examples}).

\bibliography{references_stat}{}
\bibliographystyle{plain}
\end{document}